%% file: Graph.tex
\documentclass[11pt,letterpaper]{article}
\IfFileExists{mymacros.sty}{%
  \usepackage{mymacros}%
}{%
  \input{../mymacros.sty}%
}
\usepackage{fancyhdr}
\usepackage{placeins}

\hypersetup{
  hidelinks,
  pdftitle={Traversing hyperplane arrangements with an affine objective},
  pdfauthor={Brennan Dury}
}
\fancypagestyle{plain}{
  \fancyhf{}
  \fancyfoot[C]{\thepage}

}

\renewenvironment{abstract}{%
  \small
  \begin{center}\bfseries\abstractname\end{center}%
  \noindent\ignorespaces
}{\par}
\AtBeginDocument{%
\setlength{\abovedisplayskip}{7pt plus 2pt minus 2pt}
\setlength{\belowdisplayskip}{7pt plus 2pt minus 2pt}
\setlength{\abovedisplayshortskip}{5pt plus 2pt minus 2pt}
\setlength{\belowdisplayshortskip}{5pt plus 2pt minus 2pt}
\setlength{\parskip}{2pt plus 1pt minus 1pt}
}

\newcommand{\papertitle}{Traversing hyperplane arrangements with an affine objective}
\newcommand{\paperauthor}{Brennan Dury}
\newcommand{\Ann}{\operatorname{Ann}}
\title{\papertitle}
\author{\paperauthor\\
  {\small Department of Mathematics, University of Wisconsin--Madison}\\
  {\small\href{mailto:dury@wisc.edu}{\texttt{dury@wisc.edu}}}}
\date{}

\begin{document}
\maketitle

\begin{abstract}
We give geometric constructions of directed acyclic graphs on the chambers of a locally finite affine hyperplane arrangement.
Fix a nonconstant affine objective. At each chamber where it is bounded above, a choice of flag based at the optimal face determines
\(2^r-1\) outgoing edges, where \(r\) is the codimension of that face.
The optimal value of the objective weakly increases on each step, but not always strictly. Nevertheless,
we show that every infinite directed walk in these graphs eventually reaches a chamber where the optimal value exceeds its starting value.
Additionally, we show that if an affine objective is bounded above on a closed face of a locally finite affine hyperplane arrangement, then the supremum is attained, and the optimizer set is again a closed face, generalizing the classical result for polyhedra.
\end{abstract}

\section{Introduction}
Several constructions for real hyperplane arrangements give rules for moving from one chamber to another.
In Edelman's walk, each step crosses a wall
and adds one hyperplane to the set separating the current chamber from a fixed base chamber~\cite[Proposition~1.1]{Edelman1984}.
Salvetti and Settepanella use polar coordinates to select a path of adjacent chambers from each chamber to a fixed base chamber.
Each path crosses every hyperplane at most once~\cite[Lemma~5.1]{SalvettiSettepanella2007}.
\v{C}ern\'y, Hlad\'ik, and Rada minimize a convex piecewise linear function on a finite hyperplane arrangement.
In their terminology, every chamber closure is a cell, and cells can also have lower dimension~\cite[Remark~4]{CernyHladikRada2019}.
After minimizing over a selected cell, their algorithm either stops or moves to a cell
whose minimal value is strictly smaller~\cite[\S2.3]{CernyHladikRada2019}.

We construct two directed graphs on the chambers of a locally finite affine hyperplane arrangement.
Fix a nonconstant affine objective and assign to each chamber the supremum of the objective on its closure.
At each chamber where the objective is bounded above, both constructions start a flag at the optimal face
and choose its successive faces by optimization. The walk can move in one step to any of \(2^r-1\) targets determined by the flag, where \(r\) is the codimension of the optimal face, and these chambers can be separated from the source by several hyperplanes.
The closure of each target chamber contains the optimal face, so its optimal value is at least that of the source chamber.
Notably, several successive steps can leave this value unchanged.
Nevertheless, we show that the graphs are acyclic, and for each starting chamber, every sufficiently long directed walk eventually reaches a chamber with larger optimal value than its starting chamber.
Further, if the arrangement has finitely many normal spaces, then for each starting chamber and each threshold,
every sufficiently long directed walk from that chamber ends above the threshold.

To check that the problem is well defined, we show that if an affine objective is bounded above on a closed
face of a locally finite affine hyperplane arrangement, then the supremum is attained, and the optimizer set
is again a closed face. This generalizes the usual polyhedral linear programming statement
to a class of feasible sets which need not be polyhedra. See Figure~\ref{fig:nonpolyhedral-chamber} for an example.
We also give an example in Figure~\ref{fig:locally-polyhedral-counterexample} showing that the weaker
assumption that the feasible set looks locally like a polyhedron is not sufficient. So this result shows that
if a counterexample existed, it would force infinitely many hyperplanes to concentrate on some compact set, possibly outside
the feasible set and thus surprisingly invisible to the optimization problem.

Second, the \(2^r\) chambers associated to a flag whose base has codimension \(r\) were constructed by
Aguiar and Mahajan~\cite[proof of Lemma~14.34]{AguiarMahajanArrangements}. We obtain a formula for their sign functions. The closures of these chambers contain the base, and one
of the chambers is the original chamber. So choosing a flag for each chamber defines a graph with arrows from the original chamber to the other \(2^r-1\) chambers.

Third, for any strategy assigning a flag to each chamber, suppose the chosen flags are based at the optimal faces and the graph is acyclic.
We show that for each starting chamber of finite optimal value, there is a length \(L\) such that
every directed walk from that chamber of length at least \(L\) ends in a chamber of strictly larger optimal value.
If the arrangement has finitely many normal spaces, then for each such starting chamber and each threshold,
we can choose \(L\) so that every such walk ends above the threshold.

Finally, we construct two flag selection strategies whose flags are based at the optimal face.
The Euclidean strategy uses an inner product, and the frame strategy uses covector frames, to iteratively construct
faces of the flag by optimizing over choices of the next face.
Both induce acyclic graphs, so the conclusions previously described hold.

\section{Faces and optimization in locally finite arrangements}\label{sec:arrangement-geometry}

In this paper, all affine spaces are finite dimensional.

Let $A$ be a real affine space with translation space \(V_A\), and let $\mathcal P$ be a locally finite affine hyperplane arrangement in \(A\). We write \(\Aff(A,\R)\) for the set of real affine functions on \(A\), and \(\Con(A,\R)\) for the set of real constant functions. For \(q\in\Aff(A,\R)\), define \(Z(q):=\{x\in A:q(x)=0\}\).
Then define
\(\mathfrak Q(\mathcal P):=\{q\in\Aff(A,\R)\setminus\Con(A,\R): Z(q)\in\mathcal P\}\).

For \(q\in\mathfrak Q(\mathcal P)\), let \(dq\in V_A^*\setminus\{0\}\) be its linear part.
We call \(\R dq\) the normal space of \(Z(q)\). We say that \(\mathcal P\) has finitely many
normal spaces if the set \(\{\R dq:q\in\mathfrak Q(\mathcal P)\}\) is finite.

Let
\(\Omega(\mathcal P)\) denote the chambers of $\mathcal P$.
If \(L\) is a flat, let
\(\mathcal P^L:=\{H\cap L:H\in\mathcal P,\ L\nsubseteq H,\ H\cap L\neq\emptyset\}\) be the restricted
arrangement in \(L\). We take the convention that a face of \(\mathcal P\) is a connected component of \(L\setminus|\mathcal P^L|\) for
some flat \(L\). Then if \(F\) is a face, we call \(\overline F\) a closed face.
We write \(\aff(F)\) for the affine span of a face and $\lin(F)$ for its translation space.
Write \(\Ann(F):=\{\ell\in V_A^*: \ell|_{\lin(F)}=0\}\) for its annihilator.

Define
\[
\mathcal H(F,\mathcal P):=\{H\in\mathcal P:F\subseteq H\},
\qquad
\mathfrak Q(F,\mathcal P):=\{q\in\mathfrak Q(\mathcal P):F\subseteq Z(q)\}.
\]
The map \(q\mapsto Z(q)\) takes \(\mathfrak Q(F,\mathcal P)\) onto \(\mathcal H(F,\mathcal P)\), with two
forms defining the same hyperplane exactly when one is a nonzero scalar multiple of the other. Proposition~\ref{prop:standard-face-calculus}(2) below shows that \(\mathcal H(F,\mathcal P)\) is finite, and hence that
\(\mathfrak Q(F,\mathcal P)\) is finite up to nonzero scalar rescaling.

For \(x\in A\), define the sign function
\(\sign_{\mathcal P}(x):\mathfrak Q(\mathcal P)\to\{\pm1,0\}\) by
\(\sign_{\mathcal P}(x)(q):=\sign(q(x))\).
If \(F\) is a face, by Proposition~\ref{prop:standard-face-calculus}(1), this sign function is constant on \(F\); write
\(\sign_F(q):=\sign_{\mathcal P}(x)(q)\) for any \(x\in F\).

\begin{proposition}[Face calculus]\label{prop:standard-face-calculus}
Let \(F\) be a face of a locally finite affine hyperplane arrangement \(\mathcal P\).
\begin{enumerate}
  \item Faces are exactly the nonempty fibres of $x \mapsto \sign_{\mathcal P}(x)$, and they partition \(A\). Each face \(F\) is convex and
  relatively open in \(\aff(F)\), and \(F=\relint(\overline F)\). In particular, a face is determined by its
  closure.
  \item The set \(\mathcal H(F,\mathcal P)\) is finite, and
  \[
\aff(F)=\bigcap_{H\in\mathcal H(F,\mathcal P)}H.
  \]
  Thus \(\aff(F)\) is a flat, and \(F\) is a chamber of \(\mathcal P^{\aff(F)}\).
  If \(E\) is a face with \(E\subseteq\overline F\) and \(E\neq F\), then
  \(\aff(E)\subsetneq\aff(F)\) and
  \(\mathcal H(F,\mathcal P)\subsetneq\mathcal H(E,\mathcal P)\).
  \item The closure of \(F\) is given by
  \[
  \overline F=
  \{x\in A:q(x)=0\text{ if }\sign_F(q)=0,\quad
  \sign_F(q)q(x)\geq0\text{ if }\sign_F(q)\neq0\},
  \]
  where the conditions range over \(q\in\mathfrak Q(\mathcal P)\).
  Consequently, for any face \(E\),
  \[
  E\subseteq\overline F
  \quad\Longleftrightarrow\quad
  \sign_E(q)\in\{0,\sign_F(q)\}\quad\text{for every }q\in\mathfrak Q(\mathcal P).
  \]
  \item If \(\mathfrak c\) is a chamber, the nonempty closed convex faces of \(\overline{\mathfrak c}\) are
  precisely the sets \(\overline F\) where \(F\) is a face of \(\mathcal P\) contained in \(\overline{\mathfrak c}\).
  In particular, every nonempty exposed face of \(\overline{\mathfrak c}\) is a closed face of \(\mathcal P\).
  \item For \(H_1,\ldots,H_m\in\mathcal P\), the intersection
  \(\overline F\cap H_1\cap\cdots\cap H_m\) is either empty or the closure of a unique face.
\end{enumerate}
\end{proposition}

\begin{proof}
Proposition~2.6 of \cite{DeConciniProcesiArrangements} shows that faces are convex and are precisely the
nonempty fibres of the sign function, as asserted in (1). The formula for \(\overline F\) in (3) is
established in the proof of \cite[Proposition~2.8]{DeConciniProcesiArrangements}.
Since each \(q\) has constant sign on \(E\), the same formula gives the
equivalence in (3) between \(E\subseteq\overline F\) and the stated condition on \(\sign_E\) and \(\sign_F\).

By definition, \(F\) is a chamber of the restriction to some flat \(L\). Local finiteness makes \(F\) open
in \(L\), so \(\aff(F)=L\), and convexity gives \(F=\relint(\overline F)\).
In particular, \(F\) is determined by its closure. This proves the remaining assertions in (1).

For (2), the hyperplanes containing \(F\) are exactly those containing \(L\). Since \(L\) is a flat,
their intersection is \(L=\aff(F)\). There are only finitely many such hyperplanes, since they all pass
through any chosen point of \(F\). If \(E\subseteq\overline F\), every hyperplane containing \(F\) also
contains \(E\). If moreover \(E\neq F\), then (1) implies that their sign vectors differ, and (3) gives
some \(q\) with \(\sign_E(q)=0\) and \(\sign_F(q)\neq0\). Thus \(Z(q)\) contains \(E\) but not \(F\), and
\(\aff(E)\subseteq\aff(F)\cap Z(q)\subsetneq\aff(F)\). This proves both strict inclusions in (2).

Part (4) is Remark~2.9 of \cite{DeConciniProcesiArrangements}.
For (5), put \(C:=\overline F\cap H_1\cap\cdots\cap H_m\) and suppose \(C\neq\emptyset\).
Restrict to \(L=\aff(F)\), where \(F\) is a chamber by (2). Each \(H_j\) either contains \(L\) or
cuts out a supporting hyperplane of \(\overline F\) in \(L\), by (3). Thus \(C\) is a closed convex face
of \(\overline F\). Applying (4) to \(\mathcal P^L\), we obtain \(C=\overline G\) for a face \(G\) of
\(\mathcal P^L\). By the description of faces as fibres of the sign function in (1), the faces of
\(\mathcal P^L\) are precisely the faces of \(\mathcal P\) contained in \(L\). Hence \(G\) is a face of
\(\mathcal P\), and it is uniquely determined by \(C\), by (1).
\end{proof}

We will also use the following two consequences of local finiteness.

\begin{lemma}\label{lem:local-finiteness}
Let \(\mathcal P\) be a locally finite affine hyperplane arrangement.
\begin{enumerate}
  \item If \(F\) is a face, \(x\in F\), and \(H\in\mathcal P\), then
  \(x\in H\) if and only if \(F\subseteq H\) if and only if \(\aff(F)\subseteq H\).
  There is a precompact convex open neighborhood \(U\) of \(x\) such that every arrangement hyperplane
  meeting \(U\) contains \(F\). For every
  \(q\in\mathfrak Q(\mathcal P)\setminus\mathfrak Q(F,\mathcal P)\), the sign of \(q\) on
  \(U\) is therefore \(\sign_F(q)\). If \(F\subseteq\overline{\mathfrak c}\) for a chamber
  \(\mathfrak c\), this sign is also \(\sign_{\mathfrak c}(q)\).
  \item Only finitely many closed faces meet any compact set \(B\subseteq A\). In particular, only
  finitely many chamber closures meet \(B\).
\end{enumerate}
\end{lemma}

\begin{proof}
The equivalences in (1) follow because every arrangement form has constant sign on \(F\).
Choose a ball about \(x\) with compact closure. Only finitely many arrangement hyperplanes meet it, and
shrinking the ball avoids those not passing through \(x\). This gives \(U\). Each form whose hyperplane does
not contain \(F\) has constant nonzero sign on \(U\), which must agree with \(\mathfrak c\) because
\(\mathfrak c\) meets \(U\).

For (2), choose an open ball \(U\) containing \(B\) with compact closure. Only finitely many arrangement
hyperplanes meet \(U\), and all other defining forms have fixed signs on \(U\). Thus only finitely many sign
fibres, hence faces, meet \(U\). If \(\overline F\) meets \(B\), then \(F\) meets \(U\), proving the claim.
\end{proof}

\begin{lemma}\label{lem:graded-face-intervals}
Let \(G,M\) be faces with \(G\subseteq\overline M\). Every maximal chain from \(G\) to \(M\), ordered by
closure containment, has length \(\dim M-\dim G=\codim(G)-\codim(M)\). In particular, a face \(H\) with
\(G\subsetneq\overline H\subseteq\overline M\) covers \(G\) in this interval if and only if
\(\codim(H)=\codim(G)-1\).
\end{lemma}

\begin{proof}
Choose \(x\in G\) and a neighborhood \(U\) as in Lemma~\ref{lem:local-finiteness}(1). Translate \(x\) to
the origin and keep only the finitely many hyperplanes containing \(G\). They form a finite central
arrangement whose common intersection is \(\lin(G)\), by Proposition~\ref{prop:standard-face-calculus}(2). These operations preserve partial order and dimension.

The face poset of a finite central arrangement is graded with rank given by the dimension minus the dimension
of its center~\cite[\S1.1.3]{AguiarMahajanArrangements}. That reference uses closed faces, so taking relative
interiors gives our convention. Applying this to the interval corresponding to \([G,M]\) proves the claim.
\end{proof}

We prove in Theorem~\ref{thm:locally-finite-lp-attainment} that an affine objective bounded above on a closed
face of a locally finite affine hyperplane arrangement attains its supremum, and that its optimizer set is
again a closed face.

\begin{figure}[!htbp]
\centering
\IfFileExists{Graph/figures/nonpolyhedral_chamber.tikz}{%
  \input{Graph/figures/nonpolyhedral_chamber.tikz}%
}{%
  \input{figures/nonpolyhedral_chamber.tikz}%
}
\caption{A chamber that is defined by a locally finite affine hyperplane arrangement but not a polyhedron. The boundary points
\(q_n=(n,\sqrt{n})\) escape to infinity in both \(x\) and \(y\).}
\label{fig:nonpolyhedral-chamber}
\end{figure}
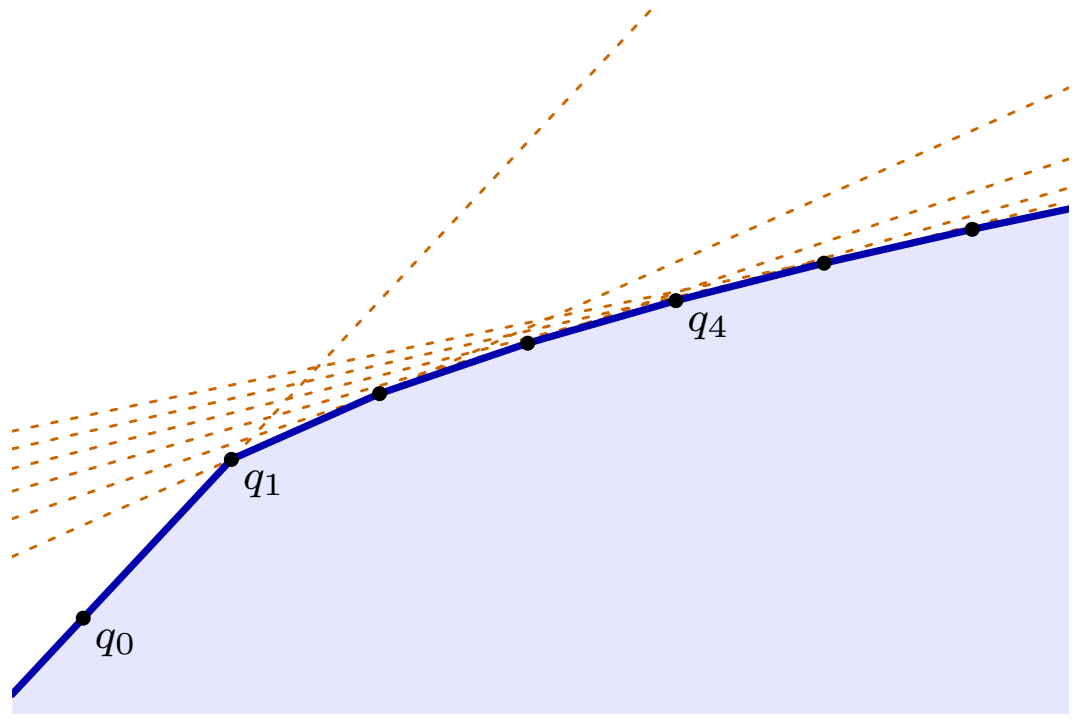

The classical attainment theorem for polyhedra does not apply directly here, since the closure of a chamber of a locally finite arrangement need not be a polyhedron.
In \(\R^2\), take the lines through consecutive points
\(q_n=(n,\sqrt n)\), \(n\geq0\), as in Figure~\ref{fig:nonpolyhedral-chamber}.
The line through \(q_n\) and \(q_{n+1}\) has slope
\(m_n=\sqrt{n+1}-\sqrt n=1/(\sqrt{n+1}+\sqrt n)\) and \(y\)-intercept \(\sqrt n-nm_n\geq\sqrt n/2\).
Since the slopes are bounded and the intercepts tend to infinity, only finitely many of these lines
meet any compact set. The points below all these lines form a chamber. Since the
slopes strictly decrease, each segment \([q_n,q_{n+1}]\) with \(n\geq1\) is an edge of its closure, whereas a
polyhedron has only finitely many edges.

\FloatBarrier

\begin{definition}\label{def:locally-polyhedral}
A set \(C\) in a real affine space \(A\) is quasipolyhedral if \(C\cap Q\) is empty or a polytope
for every polytope \(Q\subseteq A\)~\cite[p.~36, Exercise~12]{GrunbaumConvexPolytopes}.
\end{definition}

Now lower each \(q_n\) to \(p_n=(n,-2^{-n})\), as in Figure~\ref{fig:locally-polyhedral-counterexample},
and let \(H_n\) be the line through
\(p_{n}\) and \(p_{n+1}\):
\[
H_n=\left\{(x,y):y=\ell_n(x)\right\},
\qquad
\ell_n(x):=2^{-(n+1)}(x-n-2).
\]
Set
\[
C:=\{(x,y)\in\R^2:y\leq \ell_n(x)\text{ for all }n\geq0\}.
\]
Both the chamber closure in Figure~\ref{fig:nonpolyhedral-chamber} and \(C\) are quasipolyhedral.
Given a polytope \(Q\subseteq\R^2\), choose an integer \(N\geq0\) with \(x\leq N\) on \(Q\).
In each example, the lines of indices \(n\) and \(n+1\) meet at \(x=n+1\).
The second has smaller slope, so for \(n\geq N\) it lies above the first throughout \(Q\).
Intersecting \(Q\) with the first \(N+1\) closed halfplanes therefore gives its intersection with
the whole region, which is empty or a polytope.

Let \(\omega(x,y)=y\). Since \(\ell_n(x)<0\) for \(n>x\), we have \(y<0\) on \(C\).
Yet \(\omega(p_n)=-2^{-n}\to0\), so the supremum is \(0\) and is never attained.

However, the hyperplane arrangement \(\{H_n:n\geq0\}\) is not locally finite because infinitely
many hyperplanes meet the compact set \(\{(x,y) : x = 0, y \in [-1,0]\}\). 
This example shows why local finiteness is required even on compact sets outside the feasible set \(C\).

{
\setlength{\intextsep}{8pt}
\setlength{\textfloatsep}{8pt}
\setlength{\abovecaptionskip}{4pt}
\setlength{\belowcaptionskip}{0pt}
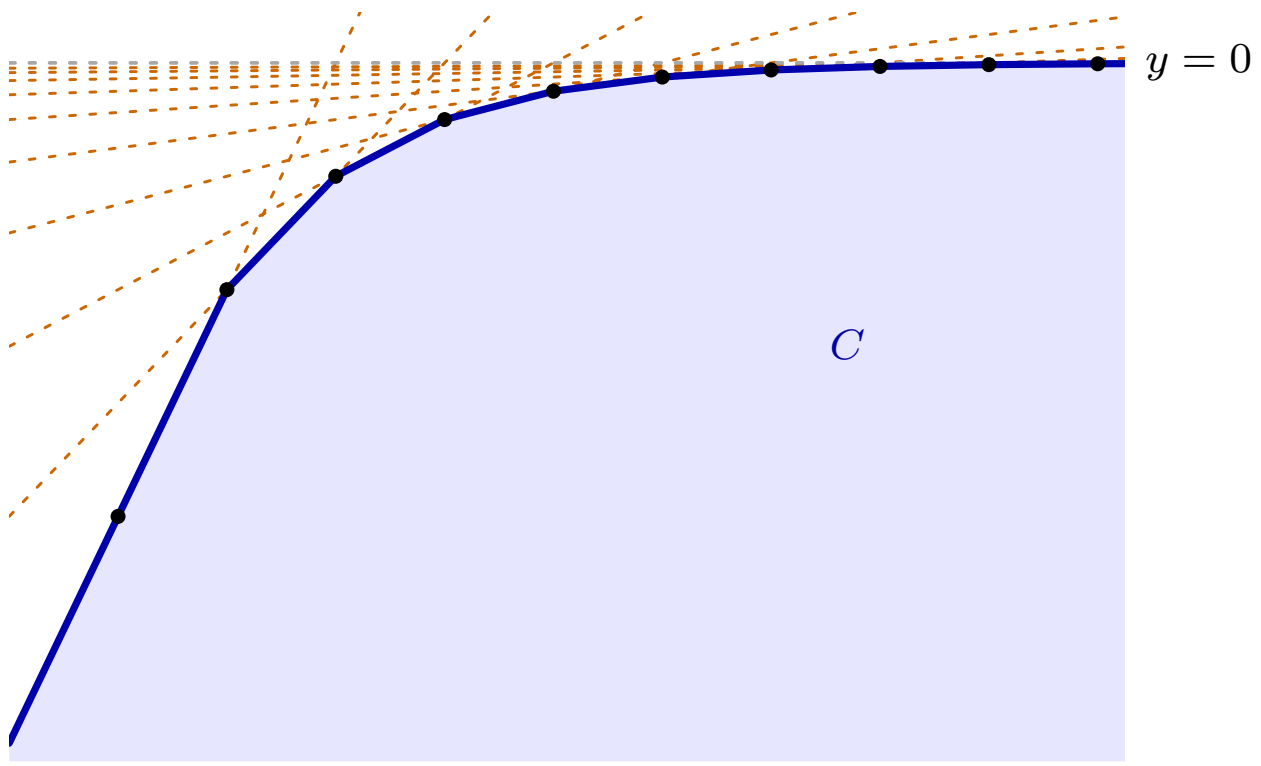
\begin{figure}[!htbp]
\centering
\scalebox{1.0}{%
  \IfFileExists{Graph/figures/unattained_optimum.tikz}{%
    \input{Graph/figures/unattained_optimum.tikz}%
  }{%
    \input{figures/unattained_optimum.tikz}%
  }%
}
\caption{A quasipolyhedral feasible set with unattained optimum.  The boundary points
\(p_n=(n,-2^{-n})\) escape to infinity in \(x\) and approach the level \(y=0\), while the supporting lines \(H_n\) accumulate
outside the feasible set \(C\).}
\label{fig:locally-polyhedral-counterexample}
\end{figure}
}
\FloatBarrier

\begin{theorem}\label{thm:locally-finite-lp-attainment}
Let \(C\) be a closed face of a locally finite affine hyperplane arrangement \(\mathcal P\).
If \(\omega\in\Aff(A,\R)\) is bounded above on \(C\), then there is
\(x\in C\) such that \(\omega(x)=\sup_{y\in C}\omega(y)\).
Moreover, the optimal set
\(\operatorname*{argmax}_{y\in C}\omega(y)\) is a closed face of \(\mathcal P\).
\end{theorem}

\begin{proof}
We argue by induction on \(\dim C\). Put \(L:=\aff(C)\). If \(\omega\) is constant on \(L\),
then every point of \(C\) is optimal. This includes the case \(\dim C=0\).
Otherwise normalize so that \(\sup_C\omega=0\), and suppose this value is not attained.
Choose \(p\in C\) and \(e\in\lin(C)\) with \(d\omega(e)=1\). Set
\(s:=p-\omega(p)e\). Then \(s\in L\) and \(\omega(s)=0\), so \(s\notin C\).

For each \(n\geq1\), choose \(x_n\in C\) with \(\omega(x_n)>-1/n\). Since \(C\) is closed and
convex, \(C\cap[x_n,s]=[x_n,z_n]\) for some \(z_n\in[x_n,s)\). The objective increases along
\([x_n,s]\), so \(\omega(z_n)\geq\omega(x_n)>-1/n\).

We claim that some \(H_n\in\mathcal P\) contains \(z_n\) and meets \([p,s]\)
(Figure~\ref{fig:planar-compact-transversal-line}). For each arrangement hyperplane through \(z_n\)
that does not contain \(L\), choose a defining affine function \(a\) with \(a\geq0\) on \(C\).
In a neighborhood of \(z_n\) in \(L\), these finitely many inequalities \(a\geq0\) define \(C\)
by Proposition~\ref{prop:standard-face-calculus}(3) and Lemma~\ref{lem:local-finiteness}(1).
If every such \(a\) satisfied \(a(s)\geq0\), then the points of \([z_n,s]\) sufficiently close
to \(z_n\) would satisfy all these inequalities and belong to \(C\). This contradicts the choice
of \(z_n\). Thus some \(a\) satisfies \(a(s)<0\). Since \(p\in C\), we have \(a(p)\geq0\),
so \(H_n:=Z(a)\) meets \([p,s]\).

By local finiteness, only finitely many arrangement hyperplanes meet the fixed segment
\([p,s]\). Hence one hyperplane \(H\) occurs for infinitely many \(n\). By
Proposition~\ref{prop:standard-face-calculus}(5), \(C\cap H\) is a closed face.
It has smaller dimension than \(C\) because \(H\) does not contain \(L\).
The corresponding \(z_n\) show that \(\sup_{C\cap H}\omega=0\). By induction this value is
attained in \(C\cap H\), contrary to our assumption.

The optimal set \(K\) is a nonempty exposed face of \(C\), hence a closed face of \(\mathcal P\)
by Proposition~\ref{prop:standard-face-calculus}(4) applied in \(L\).
\end{proof}

{
\setlength{\intextsep}{8pt}
\setlength{\textfloatsep}{8pt}
\setlength{\abovecaptionskip}{4pt}
\setlength{\belowcaptionskip}{0pt}
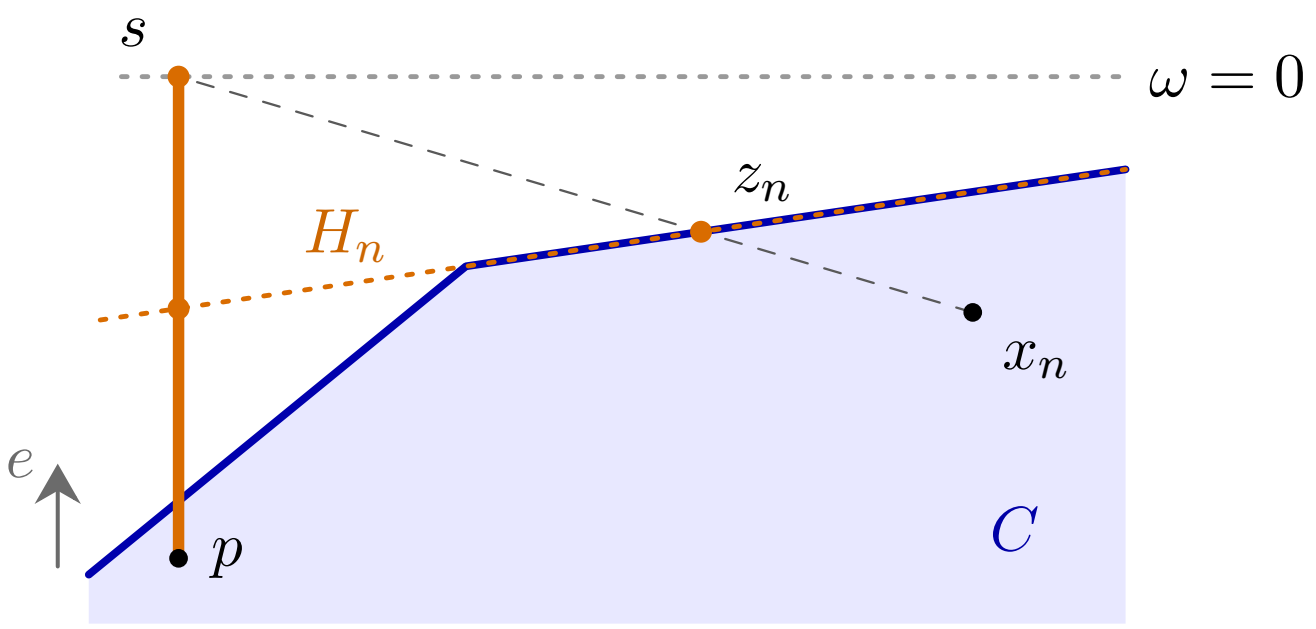
\begin{figure}[!htbp]
\centering
\scalebox{1.08}{%
  \IfFileExists{Graph/figures/attainment_fixed_segment.tikz}{%
    \input{Graph/figures/attainment_fixed_segment.tikz}%
  }{%
    \input{figures/attainment_fixed_segment.tikz}%
  }%
}
\caption{The argument in dimension two. The segment \([x_n,s]\) leaves \(C\) at \(z_n\).
The supporting line \(H_n\) through \(z_n\) meets the fixed segment \([p,s]\).}
\label{fig:planar-compact-transversal-line}
\end{figure}
}
\FloatBarrier

\section{The graph traversal}\label{sec:arrangement-flag-rules}

\subsection{The target chambers of a flag}\label{sec:flag-targets}

For a face \(F\) of \(\mathcal P\), define its chamber star by
\(\operatorname{Star}_{\Omega}(F):=\{\mathfrak c\in\Omega(\mathcal P):F\subseteq\overline{\mathfrak c}\}\).

\begin{lemma}\label{lem:finite-chamber-star}
\(\operatorname{Star}_{\Omega}(F)\) is finite.
\end{lemma}

\begin{proof}
If \(x\in F\), then every chamber in \(\operatorname{Star}_{\Omega}(F)\) has closure that meets
\(\{x\}\).  Since \(\{x\}\) is compact, only finitely many chamber closures meet \(\{x\}\) by
Lemma~\ref{lem:local-finiteness}(2).
\end{proof}

\begin{definition}\label{def:flag-of-faces}
Let \(\mathcal P\) be a locally finite affine hyperplane arrangement in \(A\).
A flag of faces is a sequence of faces \(F_r,\ldots,F_0\) such that \(F_i\subsetneq\overline{F_{i-1}}\) for all \(i\in[1,r]\), and
\(\codim(F_i)=i\).
We say that \(F_r\) is the base and \(F_0\) is the top.
\end{definition}

Let \(\mathcal P\) be a locally finite affine hyperplane arrangement in \(A\). Consider any flag of faces
with base \(F\) and top \(\mathfrak c\), written \(F=F_r,\ldots,F_0=\mathfrak c\).
For \(i\in[1,r]\), define the \(i\)-th hyperplane layer of the flag by
\(\mathcal H_i:=\mathcal H(F_i,\mathcal P)\setminus\mathcal H(F_{i-1},\mathcal P)\).

For every \(i\in[1,r]\), Proposition~\ref{prop:standard-face-calculus}(2) gives
\(\mathcal H(F_{i-1},\mathcal P)\subsetneq\mathcal H(F_i,\mathcal P)\).
Also \(\mathcal H(F_0,\mathcal P)=\emptyset\), since \(F_0=\mathfrak c\) is a chamber. Thus the
\(\mathcal H_i\) are nonempty finite sets and form a partition
\[
\mathcal H(F,\mathcal P)=\mathcal H_r\sqcup\cdots\sqcup\mathcal H_1.
\]
Equivalently, the sets \(\mathfrak Q(F_i,\mathcal P)\) form a strictly increasing chain from
\(\mathfrak Q(F_0,\mathcal P)=\emptyset\) to \(\mathfrak Q(F_r,\mathcal P)=\mathfrak Q(F,\mathcal P)\),
and \(\mathcal H_i\) records the hyperplanes of the forms added at the \(i\)-th step.

For this flag and any \(q \in \mathfrak Q(F, \mathcal P)\), define
\(\gi(q)\) as the unique index such that \(Z(q) \in \mathcal H_{\gi(q)}\).

In the finite central arrangement transverse to \(F\), the spans of the faces in the flag form
a maximal chain of flats. Aguiar and Mahajan
count the \(2^r\) flags of faces with these spans and construct them explicitly
\cite[Exercise~1.42 and proof of Lemma~14.34]{AguiarMahajanArrangements}.
Theorem~\ref{thm:opposite-flag} assigns a flag to each choice of \(r\) signs
and gives the sign of each defining form on each face of that flag.

\begin{theorem}\label{thm:opposite-flag}\prooflink{thm:opposite-flag}
Let \(\mathcal P\) be a locally finite affine hyperplane arrangement in \(A\), and let
\(\mathfrak c \in \Omega(\mathcal P)\).  Let
\(F = F_r,\ldots,F_0 = \mathfrak c\) be a flag of faces of \(\mathcal P\) with base \(F\) and top \(\mathfrak c\), with \(\codim(F)=r>0\), and denote it
by \(FL\). For each sign vector
\(\alpha=(\alpha_1,\ldots,\alpha_r)\in\{\pm 1\}^r\), there exists a chamber
\(\mathfrak b_{FL}(\alpha) \in \Omega(\mathcal P)\) and another flag of faces
\(F = F_r(\alpha),\ldots,F_0(\alpha) = \mathfrak b_{FL}(\alpha)\) with signs
\[
\sign_{F_{i}(\alpha)}(q) =
\begin{cases}
0 &: q \in \mathfrak Q(F_i, \mathcal P)\\
\alpha_{\gi(q)} \sign_{\mathfrak c}(q) &: q \in \mathfrak Q(F, \mathcal P) \setminus \mathfrak Q(F_i, \mathcal P)\\
\sign_{\mathfrak c}(q) &: q \in \mathfrak Q(\mathcal P) \setminus \mathfrak Q(F, \mathcal P)
\end{cases}.
\]
We write \(FL(\alpha)\) for this flag.  In particular, \(F_i(1^r) = F_i\) for all \(i\), so
\(\mathfrak b_{FL}(1^r)=\mathfrak c\).
\end{theorem}
\begin{proof}[\prooftarget{thm:opposite-flag}{Proof of Theorem~\ref*{thm:opposite-flag}}]
Write \(L_i=\aff(F_i)\) and set \(G_r=F_r\). We construct \(G_{r-1},\ldots,G_0\)
in that order, with \(\aff(G_i)=L_i\) at each step.

Suppose \(G_j\) has been chosen. Its span \(L_j\) is a codimension-one flat in
\(L_{j-1}\), hence a wall of the restricted arrangement. By Lemma~\ref{lem:local-finiteness}(1),
every wall meeting a sufficiently small ball about \(x\in G_j\) in \(L_{j-1}\)
contains \(\aff(G_j)=L_j\) and hence equals \(L_j\).
Let \(D_+\) be the face containing the open half-ball on the same side of \(L_j\)
as \(F_{j-1}\), and \(D_-\) the face containing the opposite half-ball.
Both faces span \(L_{j-1}\), so each has dimension \(\dim G_j+1\).
All other wall signs agree with those on \(G_j\). By parts (1) and (3) of
Proposition~\ref{prop:standard-face-calculus}, \(D_+\) and \(D_-\) are exactly
the faces with this span whose closures contain \(G_j\).
Take \(G_{j-1}=D_+\) if \(\alpha_j=1\) and \(G_{j-1}=D_-\) if \(\alpha_j=-1\). This gives the flag
\(F_i(\alpha)=G_i\), with top chamber \(G_0\).

By Proposition~\ref{prop:standard-face-calculus}(3), a nonzero sign on a face
persists on every face whose closure contains it. Thus each
\(q\in\mathfrak Q(\mathcal P)\setminus\mathfrak Q(F,\mathcal P)\)
has its sign on \(F\) throughout both flags. For \(q\in\mathfrak Q(F,\mathcal P)\),
let \(j=\gi(q)\), the unique index for which
\(q\in\mathfrak Q(F_j,\mathcal P)\setminus\mathfrak Q(F_{j-1},\mathcal P)\).
Since \(q\) is affine and vanishes on \(F_j\), it vanishes on \(L_j=\aff(F_j)\).
Thus \(q\) vanishes on \(G_i\) for \(i\ge j\), since \(G_i\subseteq L_i\subseteq L_j\).
On \(F_{j-1}\), \(q\) is nonzero. Its zero set in \(L_{j-1}\) is therefore a
hyperplane containing \(L_j\), and equals \(L_j\) because \(L_j\) has codimension one there.
The choice of side at step \(j\) gives \(\sign_{G_{j-1}}(q)=\alpha_j\sign_{F_{j-1}}(q)\).
Since \(F_{j-1}\subseteq\overline{\mathfrak c}\), Proposition~\ref{prop:standard-face-calculus}(3)
gives \(\sign_{F_{j-1}}(q)=\sign_{\mathfrak c}(q)\).
For \(i<j\), we have \(G_{j-1}\subseteq\overline{G_i}\), so the same proposition gives
\(\sign_{G_i}(q)=\sign_{G_{j-1}}(q)=\alpha_j\sign_{\mathfrak c}(q)\).
When \(\alpha=1^r\), each step selects \(F_{j-1}\) itself.
\end{proof}

\begin{definition}\label{def:flag-to-chamber}
Let \(\mathcal P\) be a locally finite affine hyperplane arrangement in \(A\). Let \(\mathfrak c \in \Omega(\mathcal P)\), and let
\(F = F_r,\ldots,F_0 = \mathfrak c\) be a flag of faces with base \(F\) and top \(\mathfrak c\), with \(r>0\). Denote the flag
by \(FL\). Theorem~\ref{thm:opposite-flag} gives a chamber
\(\mathfrak b_{FL}(\alpha)\in\operatorname{Star}_{\Omega}(F)\) for each \(\alpha\in\{\pm1\}^{r}\),
with \(\mathfrak b_{FL}(1^r)=\mathfrak c\). We call the chambers \(\mathfrak b_{FL}(\alpha)\)
with \(\alpha\ne 1^r\) the target chambers of \(FL\).
\end{definition}

The source and its \(2^r-1\) targets are mutually distinct. Indeed, if \(\alpha\neq\beta\), choose
\(i\) with \(\alpha_i\neq\beta_i\), then choose \(H\in\mathcal H_i\) and
\(q\in\mathfrak Q(F,\mathcal P)\) with \(Z(q)=H\). The sign formula in
Theorem~\ref{thm:opposite-flag} gives opposite signs for \(q\) on
\(\mathfrak b_{FL}(\alpha)\) and \(\mathfrak b_{FL}(\beta)\).

Figure~\ref{fig:two-dimensional-flag-choices} shows the target chambers in a two-dimensional example.  The
arrangement in \(A=\R^2\) has hyperplanes \(x_2=0\), \(x_1=0\), \(x_1+\sqrt{3}x_2=0\), \(x_2+x_1=0\), and
\(\sqrt{3} x_1 + x_2=0\).
For the lower left chamber, suppose the optimal face is
\(F_2=\{(0,0)\}\).
The chamber has two codimension-one faces whose closures contain \(F_2\), so there are two choices for \(F_1\).
The sign \(\alpha_2\) selects \(F_1\) or the opposite ray, and \(\alpha_1\) selects
one of the two chambers adjoining that ray.

{
\setlength{\intextsep}{6pt}
\setlength{\textfloatsep}{6pt}
\setlength{\abovecaptionskip}{3pt}
\setlength{\belowcaptionskip}{0pt}
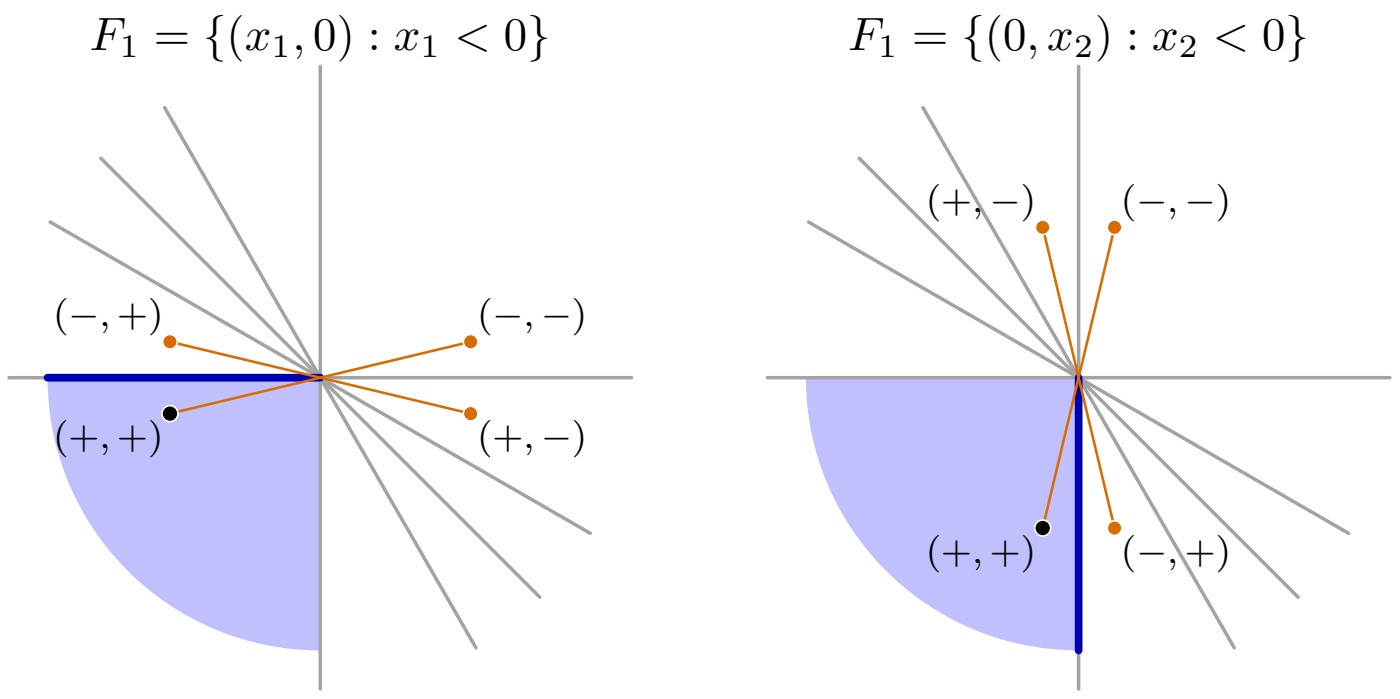
\begin{figure}[!htbp]
\centering
\scalebox{.88}{%
  \IfFileExists{figures/flag_choice_perturbations.tikz}{%
    \input{figures/flag_choice_perturbations.tikz}%
  }{%
    \input{../figures/flag_choice_perturbations.tikz}%
  }%
}
\caption{The two possible flags at the starting chamber in a two-dimensional example. In each panel
the labels give \((\alpha_1,\alpha_2)\). The black point labeled \((+,+)\) lies in the source chamber,
and the three orange points lie in its target chambers.}
\label{fig:two-dimensional-flag-choices}
\end{figure}
}
\FloatBarrier

\subsection{The directed graph of a strategy}

The target construction turns a choice of flag in each chamber into a directed graph. We require the chosen
flag to be based at the optimal face whenever the objective is bounded above. As a result, the optimal
value at the target of an edge is at least its value at the source. We then distinguish four different senses in
which the objective could be eventually increasing along a directed walk.

Let \(\mathcal P\) be a locally finite affine hyperplane arrangement in \(A\). Let
\(\omega\in\Aff(A,\R)\setminus\Con(A,\R)\). For \(\mathfrak c \in \Omega(\mathcal P)\), define the optimal value
\(\lambda(\mathfrak c, \omega) :=\sup_{x\in\overline{\mathfrak c}}\omega(x)\).

If \(\lambda(\mathfrak c,\omega)<\infty\), Theorem~\ref{thm:locally-finite-lp-attainment} shows
that the optimizer set is a closed face. Let \(F(\mathfrak c,\omega)\) be the unique face whose
closure is the optimizer set. We call \(F(\mathfrak c,\omega)\) the optimal face. A point in the interior \(\mathfrak c\) cannot be optimal
since \(d\omega\neq0\). Therefore
\[
\overline{F(\mathfrak c,\omega)}
\subseteq \partial\mathfrak c
\subseteq \overline{\mathfrak c}\cap|\mathcal P|.
\]

\begin{definition}[Flag selection strategy]\label{def:flag-strategy}
Let \(\mathcal P\) be a locally finite affine hyperplane arrangement in \(A\), and let
\(\omega\in\Aff(A,\R)\setminus\Con(A,\R)\). A flag selection strategy for \((\mathcal P,\omega)\) assigns to every chamber
\(\mathfrak c\in\Omega(\mathcal P)\) a flag with top \(\mathfrak c\),
\(\mathcal F(\mathfrak c)=\bigl(F_r(\mathfrak c),\ldots,F_0(\mathfrak c)=\mathfrak c\bigr)\)
subject to the following rules.
If \(\lambda(\mathfrak c,\omega)<\infty\), then the flag must have the optimal face as its base:
\(F_r(\mathfrak c)=F(\mathfrak c,\omega)\).
If \(\lambda(\mathfrak c,\omega)=\infty\), then \(\mathcal F(\mathfrak c)\) is the codimension-zero flag
\((F_0(\mathfrak c)=\mathfrak c)\).
\end{definition}

\begin{definition}\label{def:graph-from-flag-strategy}
Let \(\mathcal F\) be a flag selection strategy for \((\mathcal P,\omega)\). Define a directed graph \(G_{\mathcal F}\) whose
vertices are \(\Omega(\mathcal P)\) as follows. For each chamber \(\mathfrak c\), let
\(FL=\mathcal F(\mathfrak c)\), and let \(r\) be the codimension of its base face. If \(r=0\), add no edges. If
\(r>0\), add an edge from \(\mathfrak c\) to \(\mathfrak b_{FL}(\alpha)\) for each
\(\alpha\in\{\pm1\}^r\setminus\{1^r\}\).
\end{definition}
Every chamber has finite outdegree: if the base face of \(\mathcal F(\mathfrak c)\) has codimension \(r\), then there are
\(2^r-1\) outgoing edges. By Definition~\ref{def:flag-to-chamber}, these outgoing edges are distinct and are never
self-loops.

\begin{lemma}\label{lem:flag-strategy-weak-increase}
Let \(\mathcal F\) be a flag selection strategy for \((\mathcal P,\omega)\), and let
\(\mathfrak c \to \mathfrak b\) be an edge in \(G_{\mathcal F}\). Then
\(\lambda(\mathfrak b,\omega)\geq\lambda(\mathfrak c,\omega)\). If equality holds, then
\(\overline{F(\mathfrak c,\omega)}\subseteq \overline{F(\mathfrak b,\omega)}\).
\end{lemma}
\begin{proof}
Because \(\mathfrak c\) has an outgoing edge, \(\lambda(\mathfrak c,\omega)<\infty\).  By the definition of a
flag selection strategy, \(\mathcal F(\mathfrak c)\) has base \(F(\mathfrak c,\omega)\).  Every target chamber contains
that base face in its closure by Definition~\ref{def:flag-to-chamber}, so
\(F(\mathfrak c,\omega)\subseteq\overline{\mathfrak b}\).  Thus, for any \(x\in F(\mathfrak c,\omega)\),
\(\lambda(\mathfrak b,\omega)\geq \omega(x)=\lambda(\mathfrak c,\omega)\).

If equality holds, then every \(x \in F(\mathfrak c,\omega)\) satisfies \(x \in \overline{\mathfrak b}\) and
\(\omega(x) = \lambda(\mathfrak b,\omega)\). Hence
\(\overline{F(\mathfrak c,\omega)}\subseteq\overline{F(\mathfrak b,\omega)}\).
\end{proof}

\begin{definition}[Progress and escape conditions]\label{def:win-conditions}
For a flag selection strategy \(\mathcal F\), consider the following conditions on directed walks in \(G_{\mathcal F}\).
\begin{enumerate}
  \item \emph{Progress along every infinite walk.} For every chamber \(\mathfrak c_0\) with finite optimal value,
  every infinite directed walk \((\mathfrak c_0,\mathfrak c_1,\ldots)\) in \(G_{\mathcal F}\)
  eventually reaches a chamber of larger value: there is some \(i\geq0\) such that
  \(\lambda(\mathfrak c_i,\omega)>\lambda(\mathfrak c_0,\omega)\).

  \item \emph{Uniform progress.} For every chamber \(\mathfrak c_0\) with finite optimal value, there exists
  \(L_0\in\Z_{\geq1}\) such that every directed walk \((\mathfrak c_0,\ldots,\mathfrak c_L)\) in
  \(G_{\mathcal F}\) of length \(L\geq L_0\) satisfies
  \(\lambda(\mathfrak c_L,\omega)>\lambda(\mathfrak c_0,\omega)\).

  \item \emph{Escape along every infinite walk.} For every chamber \(\mathfrak c_0\) with finite optimal value, every
  infinite directed walk \((\mathfrak c_0,\mathfrak c_1,\ldots)\) in \(G_{\mathcal F}\), and every \(\eta\in\R\), there is
  some \(i\geq0\) such that \(\lambda(\mathfrak c_i,\omega)>\eta\).

  \item \emph{Uniform escape.} For every chamber \(\mathfrak c_0\) with finite optimal value and every
  \(\eta\in\R\), there exists \(L_0\in\Z_{\geq1}\) such that every directed walk
  \((\mathfrak c_0,\ldots,\mathfrak c_L)\) in \(G_{\mathcal F}\) of length \(L\geq L_0\) satisfies
  \(\lambda(\mathfrak c_L,\omega)>\eta\).
\end{enumerate}
\end{definition}

\begin{lemma}\label{lem:uniform-pathwise-progress}
Let \(\mathcal F\) be a flag selection strategy for \((\mathcal P,\omega)\). Progress along every infinite walk
is equivalent to uniform progress, and escape along every infinite walk is equivalent to uniform escape.
Either escape condition implies both progress conditions.
\end{lemma}

\begin{proof}
Each uniform condition immediately implies the corresponding condition along every infinite walk.
For the converses, fix a
chamber \(\mathfrak c_0\) with finite optimal value and a threshold \(\eta\in\R\), and assume every infinite directed walk from
\(\mathfrak c_0\) eventually has value \(>\eta\).  If \(\lambda(\mathfrak c_0,\omega)>\eta\),
Lemma~\ref{lem:flag-strategy-weak-increase} gives the uniform bound \(L_0=1\).

Assume \(\lambda(\mathfrak c_0,\omega)\leq\eta\).  If no such \(L_0\) existed, then for arbitrarily large \(L\)
there would be a directed walk \((\mathfrak c_0,\ldots,\mathfrak c_L)\) with
\(\lambda(\mathfrak c_L,\omega)\leq\eta\).  By Lemma~\ref{lem:flag-strategy-weak-increase}, every vertex of
such a directed walk has value at most \(\eta\).  The tree of finite directed walks from \(\mathfrak c_0\) whose vertices all have
value at most \(\eta\) is
finitely branching and has arbitrarily long branches, so K\"onig's infinity lemma~\cite[Lemma~8.1.2]{DiestelGraphTheory}
gives an infinite directed walk which never exceeds \(\eta\), a contradiction.

Taking \(\eta=\lambda(\mathfrak c_0,\omega)\) proves that progress along every infinite walk implies uniform progress.
Applying the argument to each \(\eta\in\R\) proves that escape along every infinite walk implies uniform escape.
Since the starting value is itself a threshold, uniform escape implies uniform progress.
\end{proof}

\begin{corollary}\label{cor:acyclic-flag-strategy-progressive}
Acyclicity of \(G_{\mathcal F}\), progress along every infinite walk, and uniform progress are equivalent.
\end{corollary}

\begin{proof}
Suppose \(G_{\mathcal F}\) is acyclic. Let
\((\mathfrak c_0,\mathfrak c_1,\ldots)\) be an infinite directed walk in \(G_{\mathcal F}\), with
\(\lambda(\mathfrak c_0,\omega)=\tau<\infty\).  Suppose, for
contradiction, that no chamber in the directed walk has \(\lambda\)-value \(>\tau\). By Lemma~\ref{lem:flag-strategy-weak-increase},
the optimal values are nondecreasing, so they all equal \(\tau\). Applying the equality case along successive edges,
we see that every chamber on the walk contains \(F(\mathfrak c_0,\omega)\) in its closure.
This face has finite star, so some chamber repeats, producing a directed cycle, contrary to acyclicity.
Thus there is progress along every infinite walk, and Lemma~\ref{lem:uniform-pathwise-progress} gives uniform progress.

Conversely, every chamber on a directed cycle has finite optimal value because it has an outgoing edge.
Following the cycle forever from a chamber of maximal value gives an infinite walk along which there is no progress.
\end{proof}

\subsection{Escape with finitely many normal spaces}
In this section, suppose \(\mathcal P\) has finitely many normal spaces.

\begin{lemma}\label{lem:finite-normal-polyhedral-chambers}\prooflink{lem:finite-normal-polyhedral-chambers}
If \(\mathcal P\) is a locally finite affine hyperplane arrangement in \(A\) with finitely many normal
spaces, then the closure of every chamber of \(\mathcal P\) is a polyhedron.
\end{lemma}

\begin{proof}[\prooftarget{lem:finite-normal-polyhedral-chambers}{Proof of Lemma~\ref*{lem:finite-normal-polyhedral-chambers}}]
For each normal space, the halfspaces in Proposition~\ref{prop:standard-face-calculus}(3) impose upper
and lower bounds on a single affine function. Taking the infimum of the upper bounds and the supremum
of the lower bounds leaves at most two inequalities per normal space, so the closure of the chamber is a polyhedron.
\end{proof}

\begin{lemma}\label{lem:shared-point-vertex}\prooflink{lem:shared-point-vertex}
Let \(\mathcal P\) be a locally finite affine hyperplane arrangement, and let
\(\mathfrak c,\mathfrak c'\in\Omega(\mathcal P)\). If
\(x\in\overline{\mathfrak c}\cap\overline{\mathfrak c'}\) is a vertex of \(\overline{\mathfrak c}\), then \(x\) is
a vertex of \(\overline{\mathfrak c'}\).
\end{lemma}

\begin{proof}[\prooftarget{lem:shared-point-vertex}{Proof of Lemma~\ref*{lem:shared-point-vertex}}]
Since \(x\) is a vertex of \(\overline{\mathfrak c}\), Proposition~\ref{prop:standard-face-calculus}(4)
implies that \(\{x\}\) is a zero-dimensional face of \(\mathcal P\). Since \(x\in\overline{\mathfrak c'}\),
the same proposition shows that \(\{x\}\) is a face of \(\overline{\mathfrak c'}\), so \(x\) is a vertex.
\end{proof}

\begin{lemma}\label{lem:finite-normal-edge-directions}
Let \(\mathcal P\) be a locally finite affine hyperplane arrangement with finitely many normal
spaces, and assume that these normal spaces span \(V_A^*\). As \(v,v'\) range over neighboring
vertices of chamber closures, the lines \(\R(v'-v)\) form a finite set.
\end{lemma}

\begin{proof}
By Proposition~\ref{prop:standard-face-calculus}(4) and (2), the direction line of each edge \(E\)
is the intersection of the spaces \(\lin(H)\) over the arrangement hyperplanes \(H\) containing \(E\).
There are only finitely many distinct spaces \(\lin(H)\), one for each normal space, so only finitely
many such lines occur.
\end{proof}

In the proof of Theorem~\ref{thm:acyclic-flag-strategy-escaping}, each target chamber contains the
source chamber's optimal face in its closure. To pass from a vertex on this face to an optimal vertex
of the target chamber, we run the simplex algorithm. To do so, we use the following standard result~\cite[Chapter~4, pp.~239--240]{LP}
which underlies the simplex algorithm.

\begin{lemma}\label{lem:better-neighbor-vertex}
Let \(C\subseteq A\) be a nonempty polyhedron, and let \(\omega\in\Aff(A,\R)\) be bounded above on \(C\).
If \(v\) is a vertex of \(C\) which is not optimal for maximizing \(\omega\), then there is a neighboring vertex
\(v'\) of \(C\) such that \(\omega(v')>\omega(v)\).
\end{lemma}

\begin{theorem}\label{thm:acyclic-flag-strategy-escaping}
Suppose \(\mathcal P\) has finitely many normal spaces, and let \(\mathcal F\) be a flag selection strategy.  If
\(G_{\mathcal F}\) is acyclic, then \(\mathcal F\) guarantees uniform escape and hence satisfies all four progress and
escape conditions.
\end{theorem}

\begin{proof}
It is enough to prove uniform escape, by Lemma~\ref{lem:uniform-pathwise-progress}.

Fix a chamber \(\mathfrak c_0\) with finite optimal value and a threshold \(\eta\in\R\).
Passing to the quotient by
\(W:=\bigcap_{q\in\mathfrak Q(\mathcal P)}\ker dq\), we may assume that the normal spaces span \(V_A^*\).
Indeed, \(\mathfrak c_0\) is invariant under translation by every vector in \(W\), so finiteness of
\(\lambda(\mathfrak c_0,\omega)\) forces \(d\omega|_W=0\). The quotient preserves
local finiteness, chambers, face incidences, and objective values, hence also the directed graph \(G_{\mathcal F}\).
Then by Lemma~\ref{lem:finite-normal-polyhedral-chambers}, every chamber closure is a polyhedron.
Any line in such a closure would be parallel to every arrangement hyperplane, so in the quotient
these polyhedra are pointed. Thus every chamber with finite optimal value has an optimal vertex in its closure~\cite[Corollary~2.11]{LP}.

By Lemma~\ref{lem:finite-normal-edge-directions}, only finitely many lines occur as directions of edges
joining two vertices of chamber closures. Fix a norm on \(V_A\) and choose unit vectors
\(u_1,\ldots,u_J\), one on each of these lines.
Let \(S:=\{|d\omega(u_j)|:j\in[1,J]\}\cap\R_{>0}\).

If \(\lambda(\mathfrak c_0,\omega)>\eta\), Lemma~\ref{lem:flag-strategy-weak-increase} gives
\(\lambda(\mathfrak c_L,\omega)>\eta\) for every directed walk from \(\mathfrak c_0\), so take \(L_0=1\).  Thus assume
\(\lambda(\mathfrak c_0,\omega)\leq\eta\).

Choose an optimal vertex \(x_0\) of \(\overline{\mathfrak c_0}\).  If \(S\neq\emptyset\), put
\(s:=\min S\).  Define
\[
B:=
\begin{cases}
\{x_0\}, & S=\emptyset,\\
\{y\in A:|y-x_0|\leq(\eta-\lambda(\mathfrak c_0,\omega))/s\}, & S\neq\emptyset.
\end{cases}
\]
The set \(B\) is compact.  By Lemma~\ref{lem:local-finiteness}(2), applied to
\(\mathcal P\), only finitely many chamber closures meet \(B\).  Let \(N\) be their number and
set \(L_0:=N+1\).

Let \((\mathfrak c_0,\ldots,\mathfrak c_L)\) be a directed walk in \(G_{\mathcal F}\) with \(L\geq L_0\).  Suppose, for a
contradiction, that \(\lambda(\mathfrak c_L,\omega)\leq\eta\).  By
Lemma~\ref{lem:flag-strategy-weak-increase}, every chamber on the directed walk has
\(\lambda(\mathfrak c_i,\omega)\leq\eta\).

Suppose \(x_i\) has been chosen as an optimal vertex of \(\overline{\mathfrak c_i}\). Since
\(\mathfrak c_{i+1}\) is a target of a flag based at \(F(\mathfrak c_i,\omega)\),
\(\overline{F(\mathfrak c_i,\omega)}\subseteq\overline{\mathfrak c_{i+1}}\).
Thus \(x_i\in\overline{\mathfrak c_{i+1}}\), and Lemma~\ref{lem:shared-point-vertex} shows that
\(x_i\) is also a vertex of \(\overline{\mathfrak c_{i+1}}\).

If \(x_i\) is already optimal on \(\overline{\mathfrak c_{i+1}}\), set
\(x_{i+1}:=x_i\), \(L_i:=0\), and \(v_{i0}:=x_i\).  Otherwise, iterating
Lemma~\ref{lem:better-neighbor-vertex} inside \(\overline{\mathfrak c_{i+1}}\) gives a walk
\[
x_i=v_{i0},\ldots,v_{iL_i}=x_{i+1}.
\]
The iteration terminates because the polyhedron has finitely many vertices and the objective strictly increases at each step.
Here \(v_{i(l-1)}\) and \(v_{il}\) are neighboring vertices of
\(\overline{\mathfrak c_{i+1}}\),
\(\omega(v_{il})>\omega(v_{i(l-1)})\) for all \(l\in[1,L_i]\), and
\(x_{i+1}\) is optimal on \(\overline{\mathfrak c_{i+1}}\).

If \(S=\emptyset\), \(\omega\) is constant along every edge joining two vertices, so no improving step
is possible. Hence \(x_j=x_0\in B\) for every \(j\).

If \(S\neq\emptyset\), each nonzero step satisfies
\(v_{il}-v_{i(l-1)}=t u_a\) for some \(t\neq0\) and some \(a\in\{1,\ldots,J\}\), and since the objective increases,
\(\sign(t)d\omega(u_a)\in S\).  Therefore
\[
\omega(x_{i+1})-\omega(x_i)
=\sum_{l=1}^{L_i}\bigl(\omega(v_{il})-\omega(v_{i(l-1)})\bigr)
\geq s\sum_{l=1}^{L_i}|v_{il}-v_{i(l-1)}|
\geq s|x_{i+1}-x_i|.
\]
For every \(j\in[0,L]\), summing these inequalities gives
\[
s\sum_{i=0}^{j-1}|x_{i+1}-x_i|
\leq \omega(x_j)-\omega(x_0)
=\lambda(\mathfrak c_j,\omega)-\lambda(\mathfrak c_0,\omega)
\leq\eta-\lambda(\mathfrak c_0,\omega).
\]
Thus \(|x_j-x_0|\leq(\eta-\lambda(\mathfrak c_0,\omega))/s\), so \(x_j\in B\) for every
\(j\).

In either case, \(B\cap\overline{\mathfrak c_j}\neq\emptyset\) for every \(j\in[0,L]\).  Hence the
chambers \(\mathfrak c_j\) lie in a fixed
finite set of size \(N\).  Since \(L\geq N+1\), some chamber repeats along the directed walk, producing a directed cycle in
\(G_{\mathcal F}\), contradicting acyclicity.  Therefore \(\lambda(\mathfrak c_L,\omega)>\eta\).  This proves
uniform escape.
\end{proof}

\section{Extending a flag one face at a time}
When the number of hyperplanes meeting at a codimension \(r\) face is greater than \(r\), there might exist no flag for which every target chamber increases the optimal value.
This allows naive strategies to form directed cycles. For the simplest example, a 2-cycle can occur. To see how, suppose that an edge
\(\mathfrak c\to\mathfrak b\) preserves the optimal face. Write \(FL=\mathcal F(\mathfrak c)\) and let
\(\alpha\) be the sign vector giving this edge. Then \(FL(\alpha)\) has the same base as \(FL\) and
has top \(\mathfrak b\), so it is an available choice at the new chamber \(\mathfrak b\).
If the strategy chooses \(FL(\alpha)\) at \(\mathfrak b\), applying \(\alpha\) again recovers every sign
because \(\alpha_i^2=1\). Theorem~\ref{thm:opposite-flag} then gives the original flag \(FL\) and the reverse edge
\(\mathfrak b\to\mathfrak c\).

In dimension two an orientation of the plane gives an acyclic strategy. If the optimal face has
codimension one, the unique flag leads to the chamber across that face, where the optimal value is larger.
If the optimal face is a vertex \(v\), draw a circle centered at \(v\) so small that it meets only
arrangement lines through \(v\). Each chamber whose closure contains \(v\) cuts out an arc.
There are two flags to choose from, one for each endpoint of the arc, as in
Figure~\ref{fig:two-dimensional-flag-choices}. Follow the arc in the chosen orientation and select
the edge of the chamber at the endpoint \(a\) reached in that direction.

The two targets with \(\alpha_2=-1\) have the antipode of \(a\) in their closures. Since \(v\) is the
unique maximizer on the source chamber closure, \(\omega(a)<\omega(v)\), so the antipode has value
greater than \(\omega(v)\).
The arcs of chambers whose optimal face is \(v\) lie in the open semicircle \(\omega<\omega(v)\), on which the chosen orientation induces a linear order.
The third target contains the next arc past \(a\) in this orientation. 
Since the orientation is chosen consistently across chambers
with the same optimal face, there cannot be a directed cycle.

Figure~\ref{fig:chamber-face-flag-example} shows a two-dimensional example with hyperplanes
\(x_2=0\), \(x_1=0\), \(x_2=-\frac{1}{\sqrt{3}}x_1\), \(x_2=-x_1\), and \(x_2=-\sqrt{3}x_1\).
The dotted arrow points in the direction of the gradient of \(\omega(x_1,x_2)=x_1+x_2\).
Chambers with optimal value \(\lambda(\mathfrak c,\omega)=0\) are
shaded light blue, while chambers with optimal value \(\lambda(\mathfrak c,\omega)>0\) are shaded light orange.
The lower left dot and darker blue shading mark the starting chamber \(\mathfrak c_0\).

Choose the clockwise orientation. The optimal face of \(\mathfrak c_0\) is the origin, and the first flag
uses the negative \(x_1\)-axis. Of its \(2^2-1=3\) target chambers, two already have optimal value greater
than zero. The third still has the origin as its optimal face, so we repeat the same move, choosing its clockwise bounding ray
\(\{(x_1,x_2):x_2=-x_1/\sqrt{3},\ x_1<0\}\). Next, one chamber with optimal value zero remains,
but now its optimal face \(\{(x_1,x_2):x_2=-x_1,\ x_1<0\}\) has codimension one. Crossing that face yields a single
new chamber, since \(2^1-1=1\), and this chamber has optimal value greater than zero. Every branch
of the traversal thus reaches a larger optimal value within three steps.

{
\setlength{\intextsep}{8pt}
\setlength{\textfloatsep}{8pt}
\setlength{\abovecaptionskip}{4pt}
\setlength{\belowcaptionskip}{0pt}
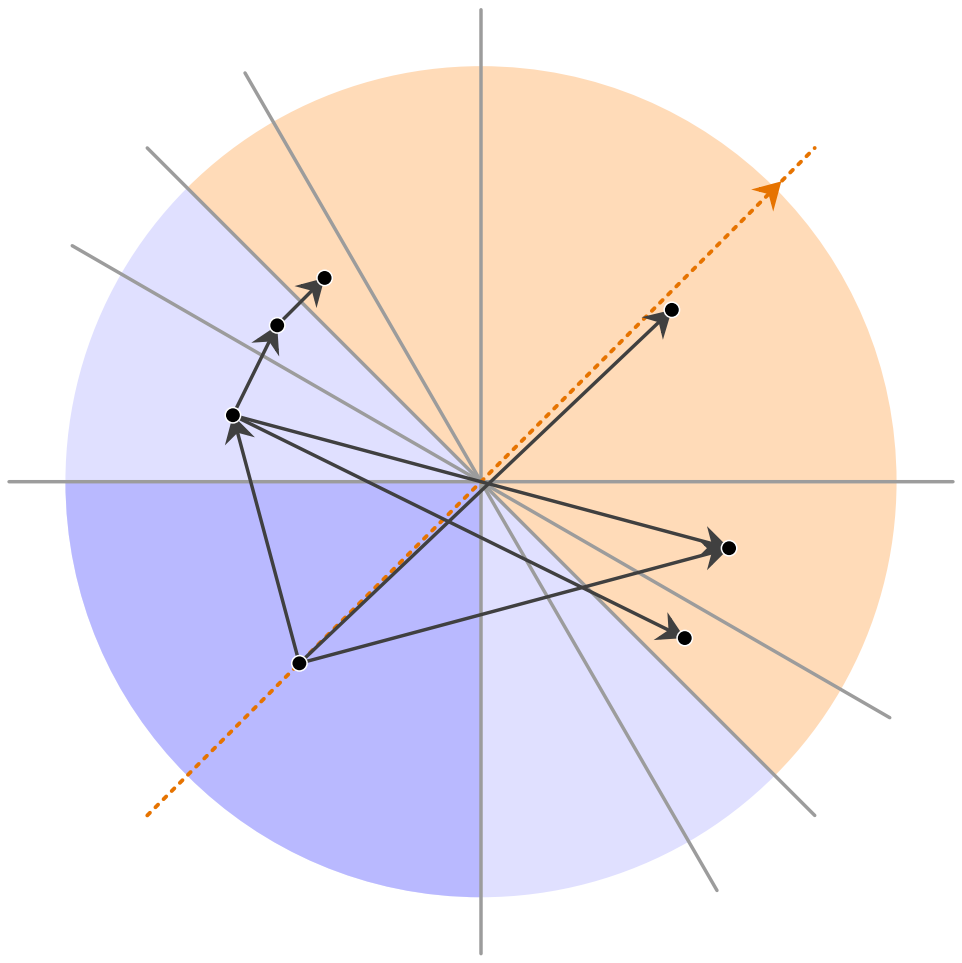
\begin{figure}[!htbp]
\centering
\scalebox{.722}{%
  \IfFileExists{figures/traversal_example.tikz}{%
    \input{figures/traversal_example.tikz}%
  }{%
    \input{../figures/traversal_example.tikz}%
  }%
}
\caption{A two-dimensional arrangement showing the directed walks from the starting chamber to chambers with larger optimal value.}
\label{fig:chamber-face-flag-example}
\end{figure}
}
\FloatBarrier

In higher dimensions an orientation of the sphere no longer singles out the next face.
The two strategies below instead choose each face of the flag iteratively, and make the choice
based on more information than an orientation.

Fix \(x\) in a face \(G\), and let \(H\) be a face with \(G\subseteq\overline H\).
We want to describe the directions \(v\in V_A\) for which \(x+tv\in H\) for all sufficiently small \(t>0\).
Whether \(v\) points into \(H\) depends only on its class in the normal quotient \(N_G:=V_A/\lin(G)\).
Write \(\pi_G:V_A\to N_G\) for the quotient map. For each \(K\in\mathcal H(G,\mathcal P)\), define
\[
K_G^\perp:=(K-x)/\lin(G)\subseteq N_G.
\]
This hyperplane is independent of \(x\). Define the transverse arrangement at \(G\) by
\[
\mathcal P_G^\perp:=\{K_G^\perp:K\in\mathcal H(G,\mathcal P)\}.
\]
Pullback by \(\pi_G\) identifies \(N_G^*\) with \(\Ann(G)\).
For \(q\in\mathfrak Q(G,\mathcal P)\), write \(dq_G^\perp\) for the induced covector on \(N_G\),
so that \(dq_G^\perp(\pi_G(v))=dq(v)\).

The directions from \(G\) into \(H\) are precisely those with \(\pi_G(v)\in R_G(H)\), where
the face \(R_G(H)\) is defined by
\[
R_G(H):=
\left\{
\xi\in N_G:
\begin{array}{ll}
dq_G^\perp(\xi)=0
& q\in\mathfrak Q(H,\mathcal P),\\
\sign_H(q)dq_G^\perp(\xi)>0
& q\in\mathfrak Q(G,\mathcal P)\setminus\mathfrak Q(H,\mathcal P)
\end{array}
\right\}.
\]

\begin{lemma}\label{lem:star-arrangement-correspondence}
The arrangement \(\mathcal P_G^\perp\) is finite, central, and essential.  The map
\[
H\longmapsto R_G(H)
\]
is an order preserving bijection from the faces \(H\) of \(\mathcal P\) satisfying
\(G\subseteq\overline H\) to the faces of \(\mathcal P_G^\perp\).  Moreover,
\[
\dim R_G(H)=\dim H-\dim G,
\qquad
\lin(R_G(H))=\lin(H)/\lin(G).
\]
\end{lemma}

\begin{proof}
By Lemma~\ref{lem:local-finiteness}(1), near any \(x\in G\) the arrangement agrees with the finite
arrangement of hyperplanes containing \(G\). After translating \(x\) to the origin, their common
intersection is \(\lin(G)\) by Proposition~\ref{prop:standard-face-calculus}(2).
The face correspondence is therefore the one in Aguiar and Mahajan~\cite[Lemma~1.38]{AguiarMahajanArrangements},
followed by quotienting by \(\lin(G)\) as in \cite[\S1.1.7]{AguiarMahajanArrangements}.
In this quotient the span of the face corresponding to \(H\) is \(\lin(H)/\lin(G)\), giving the dimension formula.
\end{proof}

In particular, extending a flag from \(G\) to a face \(H\) one dimension larger gives an open ray
\(R_G(H)\) in \(N_G\).
For \(\theta\in N_G^*\) and \(v\in V_A\), write
\(\theta(v):=\theta(\pi_G(v))\).

Let \(G\subseteq\overline{\mathfrak c}\), where \(G\) is a face and \(\mathfrak c\) is a chamber.  Define
\[
\operatorname{Cov}_G(\mathfrak c):=
\{H:G\subseteq\overline H\subseteq\overline{\mathfrak c},
\ \codim(H)=\codim(G)-1\}.
\]
By Lemma~\ref{lem:graded-face-intervals}, these are exactly the faces covering \(G\) in
\(\overline{\mathfrak c}\).

Recall that in our convention, faces are relatively open. We take extreme rays to be closed.
\begin{lemma}
\label{lem:covering-faces-generate-local-directions}
Let \(G\subseteq\overline{\mathfrak c}\).  Then \(\overline{R_G(\mathfrak c)}\) is a full-dimensional pointed polyhedral
cone in \(N_G\).  Its faces are the sets
\[
R_G(H),
\qquad
G\subseteq\overline H\subseteq\overline{\mathfrak c}.
\]
\(\overline{R_G(\mathfrak c)}\) is generated by its extreme rays \(\overline{R_G(H)}\),
one for each \(H\in\operatorname{Cov}_G(\mathfrak c)\).
\end{lemma}

\begin{proof}
Since \(R_G(\mathfrak c)\) is a chamber of the finite central arrangement \(\mathcal P_G^\perp\),
its closure is full dimensional and polyhedral. It is pointed because the arrangement is
essential: if \(\xi,-\xi\in\overline{R_G(\mathfrak c)}\), every defining form vanishes on \(\xi\), hence \(\xi=0\).
The face correspondence identifies the cone's faces with
\(R_G(H)\), and the dimension formula identifies the one-dimensional faces with the covering faces.
If \(G=\mathfrak c\), then \(N_G=\{0\}\), and the empty set of rays generates the zero cone.
Otherwise the cone is nonzero, so its extreme rays generate it~\cite[Corollary~18.5.2]{RockafellarConvexAnalysis}.
\end{proof}

Compare the source flag \(FL\) with \(FL(\alpha)\), starting at their common base.
Read \(\alpha_r,\ldots,\alpha_1\) in that order. For both strategies, the proof of acyclicity will focus on
the first \(-1\). At each \(+1\) before the first \(-1\), the two
flags choose the same next face. Then at the first \(-1\), their next faces determine opposite rays
in the normal quotient of their last common face.

\begin{lemma}\label{lem:first-flag-change}\prooflink{lem:first-flag-change}
Assume the setting of Theorem~\ref{thm:opposite-flag}. Suppose \(\alpha \neq 1^r\), and let \(i^*\) be the
largest index such that \(\alpha_{i^*}=-1\).  Then \(F_i=F_i(\alpha)\) for \(i^*\leq i\leq r\).
In particular,
\(F_i\subseteq\overline{\mathfrak b_{FL}(\alpha)}\cap\overline{\mathfrak c}\) for \(i^*\leq i\leq r\).

The two rays satisfy
\[
R_{F_{i^*}}(F_{i^*-1}(\alpha))=-R_{F_{i^*}}(F_{i^*-1}).
\]
\end{lemma}
\begin{proof}[\prooftarget{lem:first-flag-change}{Proof of Lemma~\ref*{lem:first-flag-change}}]
Fix \(i\in[i^*,r]\).  We compare the sign descriptions in Theorem~\ref{thm:opposite-flag}.  If
\(q\in\mathfrak Q(F_i,\mathcal P)\), then both faces have sign \(0\). If
\(q\in\mathfrak Q(F,\mathcal P)\setminus\mathfrak Q(F_i,\mathcal P)\), then
\(Z(q)\in\mathcal H_k\) for some \(k>i\), so \(\gi(q)=k>i\geq i^*\).  Hence
\(\alpha_{\gi(q)}=1\) by maximality of \(i^*\).  Thus the two nonzero signs are both
\(\sign_{\mathfrak c}(q)\).  For
\(q\notin\mathfrak Q(F,\mathcal P)\), both signs are \(\sign_{\mathfrak c}(q)\).  Hence
\(F_i\) and \(F_i(\alpha)\) have the same sign pattern, so \(F_i=F_i(\alpha)\).

The inclusion follows because \(F_i(\alpha) \subseteq \overline{F_0(\alpha)}
=\overline{\mathfrak b_{FL}(\alpha)}\) and
\(F_i\subseteq \overline{F_0}=\overline{\mathfrak c}\).

For \(q\in\mathfrak Q(F_{i^*},\mathcal P)\setminus\mathfrak Q(F_{i^*-1},\mathcal P)\), one has
\(\gi(q)=i^*\), so the signs of \(F_{i^*-1}\) and \(F_{i^*-1}(\alpha)\) on \(q\) are opposite.
The same forms vanish on both faces. Consequently
\[
R_{F_{i^*}}(F_{i^*-1}(\alpha))=-R_{F_{i^*}}(F_{i^*-1}).
\]
\end{proof}

\section{The Euclidean strategy}\label{sec:euclidean-strategy}

In this section, fix an inner product on \(V_A\). For convenience, write \(\ell=d\omega\in V_A^*\).

Return to the starting chamber \(\mathfrak c_0\) of Figure~\ref{fig:chamber-face-flag-example}, with the
standard inner product. The two rays bounding \(\mathfrak c_0\) meet the unit circle at \((-1,0)\) and
\((0,-1)\), as shown in Figure~\ref{fig:euclidean-unit-circle}. On the southwest arc,
\(\ell(x_1,x_2)=x_1+x_2\) attains its maximum of \(-1\) at both endpoints, each corresponding to an extreme ray.
When the endpoint values differ, we choose the ray through the endpoint with larger value.
Here they agree, so we use a generic second covector to break the tie, consistently chosen for all
chambers with the same optimal face. Here \(\ell'_F(x_1,x_2)=x_2\) works. It takes the value \(0\)
at \((-1,0)\) and \(-1\) at \((0,-1)\), so we choose the negative \(x_1\)-axis.

In the neighboring chamber across this ray, the origin is still optimal. The two bounding rays
have unit generators \((-1,0)\) and \((-\sqrt{3}/2,1/2)\). The latter gives
\(\ell(-\sqrt{3}/2,1/2)=(1-\sqrt{3})/2>-1\), so maximizing \(\ell\) alone selects the next
clockwise ray in Figure~\ref{fig:chamber-face-flag-example}.

{
\setlength{\intextsep}{8pt}
\setlength{\textfloatsep}{8pt}
\setlength{\abovecaptionskip}{4pt}
\setlength{\belowcaptionskip}{0pt}
\begin{figure}[!htbp]
\centering
\IfFileExists{Graph/EuclideanStrategy/figures/unit_circle_selection.tikz}{%
  \input{Graph/EuclideanStrategy/figures/unit_circle_selection.tikz}%
}{%
  \IfFileExists{EuclideanStrategy/figures/unit_circle_selection.tikz}{%
    \input{EuclideanStrategy/figures/unit_circle_selection.tikz}%
  }{%
    \input{figures/unit_circle_selection.tikz}%
  }%
}
\caption{Both endpoints of the southwest arc lie on the level set
\(\ell=-1\). The secondary covector at the origin selects \((-1,0)\) by maximizing \(x_2\).}
\label{fig:euclidean-unit-circle}
\end{figure}
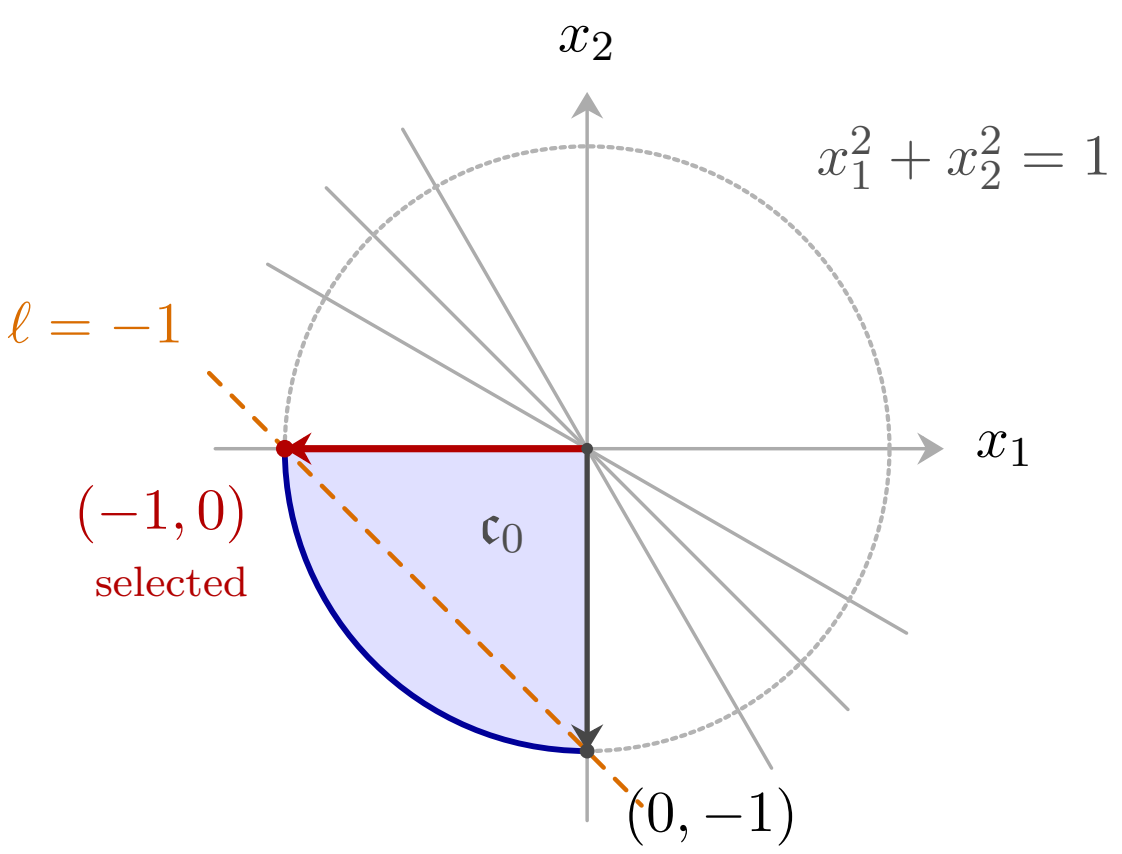
}
\FloatBarrier

\subsection{Flag selection}

We proceed by generalizing this construction to higher dimensions.
For a face \(F\) of positive codimension, identify
\(N_F=V_A/\lin(F)\) with \(\lin(F)^\perp\) using the inner product.

\begin{definition}\label{def:euclidean-generator}
Let \(H\) be a face with \(F\subseteq\overline H\) and \(\dim H=\dim F+1\).
Define the generator of \(H\) over \(F\) to be the unique unit vector
\[
u(H,F)\in\lin(H)\cap\lin(F)^\perp
\]
whose image in \(N_F\) belongs to the ray \(R_F(H)\) of
Lemma~\ref{lem:star-arrangement-correspondence}.
\end{definition}

Equivalently, for \(x\in F\), one has \(x+t\,u(H,F)\in H\) for all sufficiently small
\(t>0\).

For \(F\subseteq\overline{\mathfrak c}\), with \(F\neq\mathfrak c\), the available covering faces
are \(\operatorname{Cov}_F(\mathfrak c)\). This set is nonempty and finite by
Lemma~\ref{lem:covering-faces-generate-local-directions}.
For \(\eta\in V_A^*\), define
\[
P(F,\mathfrak c,\eta)
:=\max_{H\in\operatorname{Cov}_F(\mathfrak c)}\eta\bigl(u(H,F)\bigr).
\]
The maximizing covering face for \(\eta=\ell\) need not be unique. For each optimal face \(F\), we choose
a second covector \(\ell'_F\in N_F^*\) to break ties throughout every flag based at \(F\).
At a subsequent face \(G\) with \(F\subseteq\overline G\), the unit generators belong to
\(\lin(G)^\perp\subseteq\lin(F)^\perp\simeq N_F\), so \(\ell'_F\) can evaluate them as well.

Call a covector admissible at \(F\) if, for every chamber \(\mathfrak c\) with optimal face \(F\)
and every face \(G\ne\mathfrak c\) with \(F\subseteq\overline G\subseteq\overline{\mathfrak c}\),
it takes pairwise distinct values on the generators \(u(H,G)\), \(H\in\operatorname{Cov}_G(\mathfrak c)\).

\begin{lemma}\label{lem:euclidean-generic-tiebreaker}
For each optimal face \(F\), the admissible covectors form a nonempty Zariski open subset of \(N_F^*\).
\end{lemma}

\begin{proof}
There are finitely many faces whose closures contain \(F\), by the face correspondence
of Lemma~\ref{lem:star-arrangement-correspondence}. For each \(\mathfrak c,G\) in the definition
and each pair of distinct covers \(H,H'\), the tie
\[
\eta\bigl(u(H,G)-u(H',G)\bigr)=0
\]
defines a proper hyperplane in \(N_F^*\), since the two unit generators are distinct.
The admissible covectors are exactly the complement of these finitely many hyperplanes.
\end{proof}

Choose an admissible covector \(\ell'_F\) independently of the chamber at each optimal face \(F\).

\begin{definition}\label{def:euclidean-selection}
Let \(F\) be the optimal face of \(\mathfrak c\). For \(F\subseteq\overline G\subseteq\overline{\mathfrak c}\)
with \(G\ne\mathfrak c\), let \(Q(G,\mathfrak c)\) be the unique \(H\in\operatorname{Cov}_G(\mathfrak c)\) maximizing
\[
\bigl(\ell(u(H,G)),\,\ell'_F(u(H,G))\bigr)
\]
in lexicographic order. Denote the maximizing pair by \(S(G,\mathfrak c)\).
Its first coordinate is \(P(G,\mathfrak c,\ell)\).
\end{definition}

\begin{definition}\label{def:euclidean-flag}
Let \(\lambda(\mathfrak c,\omega)<\infty\), and put
\(r=\codim F(\mathfrak c,\omega)\). Define recursively
\[
F_r^{\mathrm E}(\mathfrak c):=F(\mathfrak c,\omega),
\qquad
F_{i-1}^{\mathrm E}(\mathfrak c)
:=
Q(F_i^{\mathrm E}(\mathfrak c),\mathfrak c)
\quad (i=r,\ldots,1).
\]
Denote the resulting flag by
\[
FL_{\mathrm E}(\mathfrak c)
:=\bigl(F_r^{\mathrm E}(\mathfrak c),\ldots,F_0^{\mathrm E}(\mathfrak c)=\mathfrak c\bigr).
\]
If \(\lambda(\mathfrak c,\omega)=\infty\), choose the codimension-zero flag
\((\mathfrak c)\). Denote this flag selection strategy by \(\mathcal F_{\mathrm E}\), and let
\(G_{\mathrm E}:=G_{\mathcal F_{\mathrm E}}\) be its graph as in
Definition~\ref{def:graph-from-flag-strategy}.
\end{definition}

Also for a chamber with finite optimal value, put
\[
u_i(\mathfrak c)
:=
u(F_{i-1}^{\mathrm E}(\mathfrak c),F_i^{\mathrm E}(\mathfrak c)).
\]
These vectors form an orthonormal basis of
\(\lin(F(\mathfrak c,\omega))^\perp\) such that each vector points from the current face
into the next face of the selected flag.

\subsection{Strict increase of the potential}

Now we have constructed the strategy, so next we show that it is acyclic. For a chamber with finite optimal value and optimal face \(F\), write
\[
S_i^{\mathrm E}(\mathfrak c)
:=
S(F_i^{\mathrm E}(\mathfrak c),\mathfrak c)
=\bigl(\ell(u_i(\mathfrak c)),\,\ell'_F(u_i(\mathfrak c))\bigr).
\]

\begin{definition}\label{def:euclidean-potential}
If \(\lambda(\mathfrak c,\omega)<\infty\), put \(r=\codim F(\mathfrak c,\omega)\) and define the potential by
\[
P_{\mathrm E}^*(\mathfrak c)
:=
\bigl(\lambda(\mathfrak c,\omega),\,
      \dim F(\mathfrak c,\omega),\,
      S_r^{\mathrm E}(\mathfrak c),\ldots,S_2^{\mathrm E}(\mathfrak c)\bigr).
\]
The last face of the flag is forced to be \(\mathfrak c\), so the sequence stops at codimension two.

Compare these tuples lexicographically.
If their first two coordinates agree, their optimal
faces have the same codimension, so the remaining tuples have the same length.

If \(\lambda(\mathfrak c,\omega)=\infty\), set \(P_{\mathrm E}^*(\mathfrak c)=+\infty\),
ordered above every chamber with finite optimal value.
\end{definition}

Observe that \(\omega\) decreases in every direction from the optimal face into any covering face; otherwise, that covering face must also be optimal. So \(P(F_r^{\mathrm E}(\mathfrak c),\mathfrak c,\ell)<0\). We solve a related optimization problem at each step, and our next goal is to show that similarly, \(P(F_i^{\mathrm E}(\mathfrak c),\mathfrak c,\ell)<0\) for all \(1 \leq i \leq r\). We can then show that the potential increases strictly
along every directed edge, which rules out directed cycles.

\begin{lemma}\label{lem:euclidean-two-intermediate-faces}
Suppose \(F\subseteq\overline K\subseteq\overline{\mathfrak c}\) and
\(\dim K=\dim F+2\). There are exactly two faces \(H,H'\) such that
\[
F\subseteq\overline H\subseteq\overline K,\qquad \dim H=\dim F+1.
\]
Let
\[
u=u(H,F),\qquad v=u(K,H),\qquad w=u(H',F).
\]
Then
\[
w=au+bv,\qquad b>0,\qquad a^2+b^2=1,\qquad a<1.
\]
\end{lemma}

\begin{proof}
By Lemma~\ref{lem:covering-faces-generate-local-directions},
\(\overline{R_F(K)}\) is a two-dimensional pointed polyhedral cone, and its two extreme rays
correspond precisely to \(H,H'\). Under the identification \(N_F\simeq\lin(F)^\perp\), the cone spans
\(\lin(K)\cap\lin(F)^\perp\). The vectors \(u,v\) form an orthonormal basis of this plane,
so \(w=au+bv\) for some \(a,b \in \R\).

Since \(\overline{R_F(H)}\) is a boundary face of the convex cone \(\overline{R_F(K)}\),
its span \(\lin(H)/\lin(F)=\R u\) is a supporting line. The generator \(v\) points into
the halfplane containing the cone. The other boundary ray lies strictly in that halfplane:
it is distinct from \(\R_{\geq0}u\), and pointedness excludes the opposite ray.
Thus \(b=\langle w,v\rangle>0\). Since \(w\) is a unit vector,
\(a^2+b^2=1\), which also gives \(a<1\).
\end{proof}

{
\setlength{\intextsep}{8pt}
\setlength{\textfloatsep}{8pt}
\setlength{\abovecaptionskip}{4pt}
\setlength{\belowcaptionskip}{0pt}
\begin{figure}[!htbp]
\centering
\IfFileExists{Graph/EuclideanStrategy/figures/two_intermediate_faces.tikz}{%
  \input{Graph/EuclideanStrategy/figures/two_intermediate_faces.tikz}%
}{%
  \IfFileExists{EuclideanStrategy/figures/two_intermediate_faces.tikz}{%
    \input{EuclideanStrategy/figures/two_intermediate_faces.tikz}%
  }{%
    \input{figures/two_intermediate_faces.tikz}%
  }%
}
\caption{The cone in Lemma~\ref{lem:euclidean-two-intermediate-faces},
with boundary rays \(R_F(H)\) and \(R_F(H')\).
The dashed segment from \(a u\) to \(w\) is \(b v\), so \(w=a u+b v\).
The coefficient \(a\) is positive on the left and negative on the right, but in both cases \(b>0\).
}
\label{fig:euclidean-two-intermediate-faces}
\end{figure}
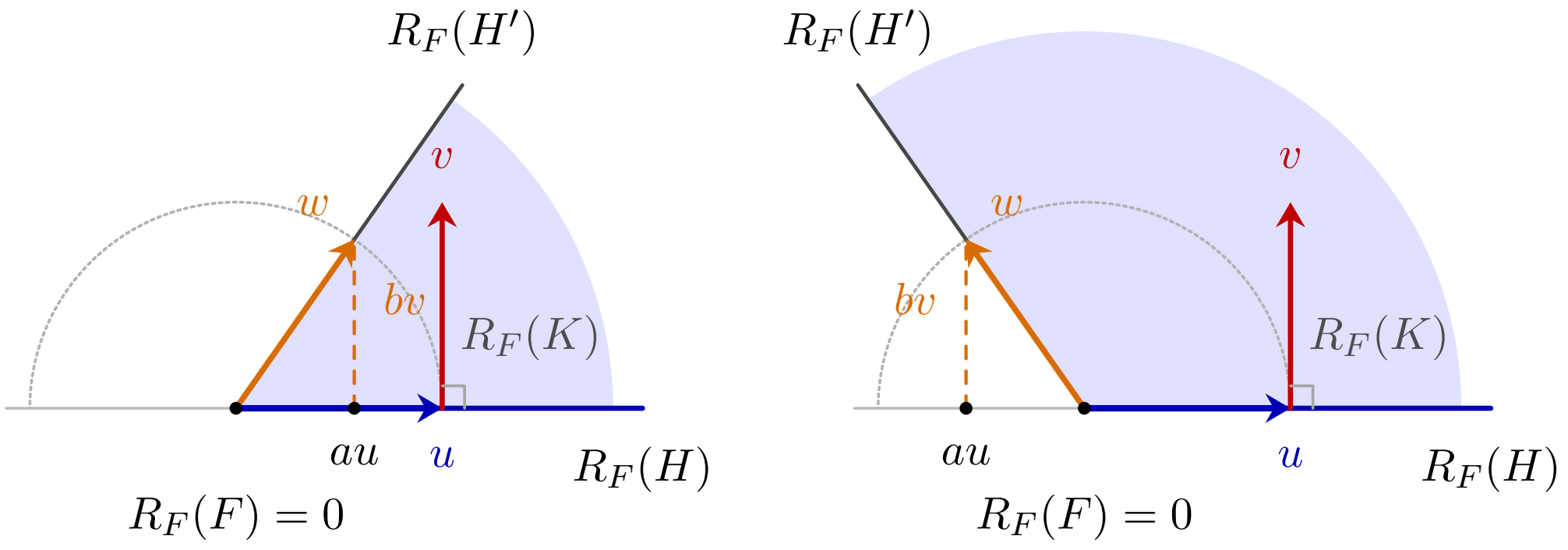
}
\FloatBarrier

\begin{lemma}\label{lem:euclidean-negative-derivatives}
For every chamber \(\mathfrak c\) with \(\lambda(\mathfrak c,\omega)<\infty\),
put \(r=\codim F(\mathfrak c,\omega)\). Then
\[
\ell(u_i(\mathfrak c))
=P(F_i^{\mathrm E}(\mathfrak c),\mathfrak c,\ell)<0
\qquad (1\leq i\leq r).
\]
\end{lemma}

\begin{proof}
Write \(F_i=F_i^{\mathrm E}(\mathfrak c)\) and \(u_i=u_i(\mathfrak c)\).
For the base case \(i=r\), choose \(x\in F_r=F(\mathfrak c,\omega)\).
For sufficiently small \(t>0\), the point \(x+tu_r\) belongs to \(F_{r-1}\), hence to
\(\overline{\mathfrak c}\setminus\overline{F_r}\). Therefore
\[
\omega(x+tu_r)<\omega(x)=\lambda(\mathfrak c,\omega),
\]
so \(\ell(u_r)<0\).

Suppose \(2\leq i\leq r\) and \(\ell(u_i)<0\). Apply
Lemma~\ref{lem:euclidean-two-intermediate-faces} to
\[
F_i\subseteq\overline{F_{i-1}}\subseteq\overline{F_{i-2}}.
\]
Let \(H'\neq F_{i-1}\) be the other intermediate face, and let \(w=u(H',F_i)\).
Then
\[
w=a u_i+b u_{i-1},\qquad b>0,\quad a<1.
\]
Since \(H'\in\operatorname{Cov}_{F_i}(\mathfrak c)\), it was an available candidate when
\(F_{i-1}\) was selected. So the selection rule gives
\(\ell(w)\leq\ell(u_i)\). Consequently,
\[
a\ell(u_i)+b\ell(u_{i-1})\leq\ell(u_i),
\qquad
\ell(u_{i-1})\leq\frac{1-a}{b}\ell(u_i)<0.
\]
Descending induction proves the assertion.
\end{proof}

\begin{proposition}\label{prop:euclidean-potential-increase}
For every edge \(\mathfrak c\to\mathfrak b\) in \(G_{\mathrm E}\),
\[
P_{\mathrm E}^*(\mathfrak b)>_{\mathrm{lex}}P_{\mathrm E}^*(\mathfrak c).
\]
\end{proposition}

\begin{proof}
The source has finite optimal value, since chambers with infinite optimal value have no outgoing edges.
By Lemma~\ref{lem:flag-strategy-weak-increase},
\(\lambda(\mathfrak b,\omega)\geq\lambda(\mathfrak c,\omega)\).
If the target has infinite optimal value, its potential is \(+\infty\), so the conclusion holds.
If both optimal values are finite and the inequality is strict, the first coordinate of the potential increases.
Otherwise the same lemma gives
\[
\overline{F(\mathfrak c,\omega)}
\subseteq\overline{F(\mathfrak b,\omega)}.
\]
If the optimal faces have different dimensions, the second coordinate increases.
We may therefore assume they are the same face \(F\), of codimension \(r\).

Write \(F_i=F_i^{\mathrm E}(\mathfrak c)\) and \(G_i=F_i^{\mathrm E}(\mathfrak b)\).
There is \(\alpha\in\{\pm1\}^r\setminus\{1^r\}\) such that
\(\mathfrak b=\mathfrak b_{FL_{\mathrm E}(\mathfrak c)}(\alpha)\).
Besides the two selected flags \((F_r,\ldots,F_0)\) and \((G_r,\ldots,G_0)\), consider the
flag \((F_r(\alpha),\ldots,F_0(\alpha))\) given by
Theorem~\ref{thm:opposite-flag}. Put
\[
i^*:=\max\{i:\alpha_i=-1\}.
\]
Lemma~\ref{lem:first-flag-change} gives
\[
F_i(\alpha)=F_i
\quad(i^*\leq i\leq r),
\qquad
u(F_{i^*-1}(\alpha),F_{i^*})=-u_{i^*}(\mathfrak c).
\]
Indeed,
\[
R_{F_{i^*}}\!\left(F_{i^*-1}\right)
\quad\text{and}\quad
R_{F_{i^*}}\!\left(F_{i^*-1}(\alpha)\right)
\]
are opposite, so their unit generators are negatives of one another.

We compare the two selected flags starting from their common base \(F_r=G_r=F\), so both
chambers use the same secondary covector \(\ell'_F\).
Suppose first that they make different choices at some level \(j>i^*\), and take the first
such level in the descending recursion. Thus
\[
F_a=G_a\quad(j\leq a\leq r),
\qquad F_{j-1}\neq G_{j-1}.
\]
Since \(j-1\geq i^*\), Lemma~\ref{lem:first-flag-change} gives
\(F_{j-1}=F_{j-1}(\alpha)\subseteq\overline{\mathfrak b}\).
It is therefore an available covering face at \(G_j=F_j\)
in the target's selection. Its pair of covector values is \(S_j^{\mathrm E}(\mathfrak c)\).
The target selects a different face, its unique lexicographic maximizer, so
\[
S_j^{\mathrm E}(\mathfrak b)>_{\mathrm{lex}}S_j^{\mathrm E}(\mathfrak c).
\]
All earlier pairs \(S_a^{\mathrm E}\), \(a>j\), agree. Hence the target potential is strictly larger.

It remains to consider the case in which the selected flags agree through \(F_{i^*}\):
\[
G_a=F_a\qquad(i^*\leq a\leq r).
\]
Now \(F_{i^*-1}(\alpha)\) is an available covering face of \(G_{i^*}\) in
\(\overline{\mathfrak b}\), with generator \(-u_{i^*}(\mathfrak c)\).
The target maximizes the first coordinate before breaking ties, so
\[
\ell(u_{i^*}(\mathfrak b))
\geq-\ell(u_{i^*}(\mathfrak c))
>0,
\]
where the strict inequality follows from
Lemma~\ref{lem:euclidean-negative-derivatives} applied to \(\mathfrak c\).
This contradicts the same lemma applied to \(\mathfrak b\), so this case cannot occur.
\end{proof}

\begin{theorem}\label{thm:euclidean-progress-and-escape}\prooflink{thm:euclidean-progress-and-escape}
Let \(\mathcal P\) be a locally finite affine hyperplane arrangement in \(A\), and let
\(\omega\colon A\to\R\) be a nonconstant affine objective. Fix an inner product on \(V_A\).
At each optimal face \(F\), there exists an admissible covector \(\ell'_F\).
For every such family \((\ell'_F)_F\), the graph \(G_{\mathrm E}\) is acyclic and the strategy
\(\mathcal F_{\mathrm E}\) guarantees both progress along every infinite walk and uniform progress. If \(\mathcal P\) has finitely many
normal spaces, then \(\mathcal F_{\mathrm E}\) also guarantees both escape along every infinite walk and uniform escape.
\end{theorem}

\begin{proof}[\prooftarget{thm:euclidean-progress-and-escape}{Proof of Theorem~\ref*{thm:euclidean-progress-and-escape}}]
Lemma~\ref{lem:euclidean-generic-tiebreaker} gives an admissible covector at each optimal face. For any such choice,
Proposition~\ref{prop:euclidean-potential-increase} makes the potential strictly increase
along every directed edge, which excludes directed cycles.
The progress and escape assertions now follow from
Corollary~\ref{cor:acyclic-flag-strategy-progressive} and
Theorem~\ref{thm:acyclic-flag-strategy-escaping}.
\end{proof}

\section{The frame strategy}\label{sec:our-strategy}
The Euclidean strategy uses an inner product and breaks ties with a covector chosen at each optimal face.
We provide an alternative which instead depends on a choice of frame of covectors.

Return to the chamber \(\mathfrak c_0\) in Figure~\ref{fig:chamber-face-flag-example}. We will produce the same choice of flags with a different perspective.
Since \(\omega\) is negative on \(\overline{\mathfrak c_0}\setminus\{0\}\), the line
\(\omega=-1\) cuts its two bounding rays at the endpoints of a segment.
The covector \(\omega'(x_1,x_2)=x_2\) takes its larger value at the endpoint \(p\) on the
negative \(x_1\)-axis, so the flag uses the ray \(F_1=\R_{>0}p\) (Figure~\ref{fig:two-dimensional-score}).

We can see on the same line what happens in the neighboring chamber across \(F_1\).
Write \(p+v\) for the other endpoint of the first segment. Since \(p\) maximizes \(\omega'\),
we have \(\omega'(v)<0\). The neighboring segment runs from \(p\) in the direction \(-v\)
to its other endpoint \(p-tv\), where \(t>0\). Thus \(\omega'(p-tv)>\omega'(p)\). So \(\omega'\) increases in the next chamber.

{
\setlength{\intextsep}{8pt}
\setlength{\textfloatsep}{8pt}
\setlength{\abovecaptionskip}{4pt}
\setlength{\belowcaptionskip}{0pt}
\begin{figure}[!htbp]
\centering
\scalebox{.95}{%
  \IfFileExists{figures/two_dimensional_score.tikz}{%
    \input{figures/two_dimensional_score.tikz}%
  }{%
    \input{../figures/two_dimensional_score.tikz}%
  }%
}
\caption{
The shaded lower left chamber is \(\mathfrak c_0\). The orange segment is
\(\overline{\mathfrak c_0}\cap\{\omega=-1\}\), with endpoints \(p\) and \(p+v\). The red arrow
starts at \(p\) on the common boundary ray \(F_1\) and ends at the maximizing point \(p-tv\)
of the neighboring slice.}
\label{fig:two-dimensional-score}
\end{figure}
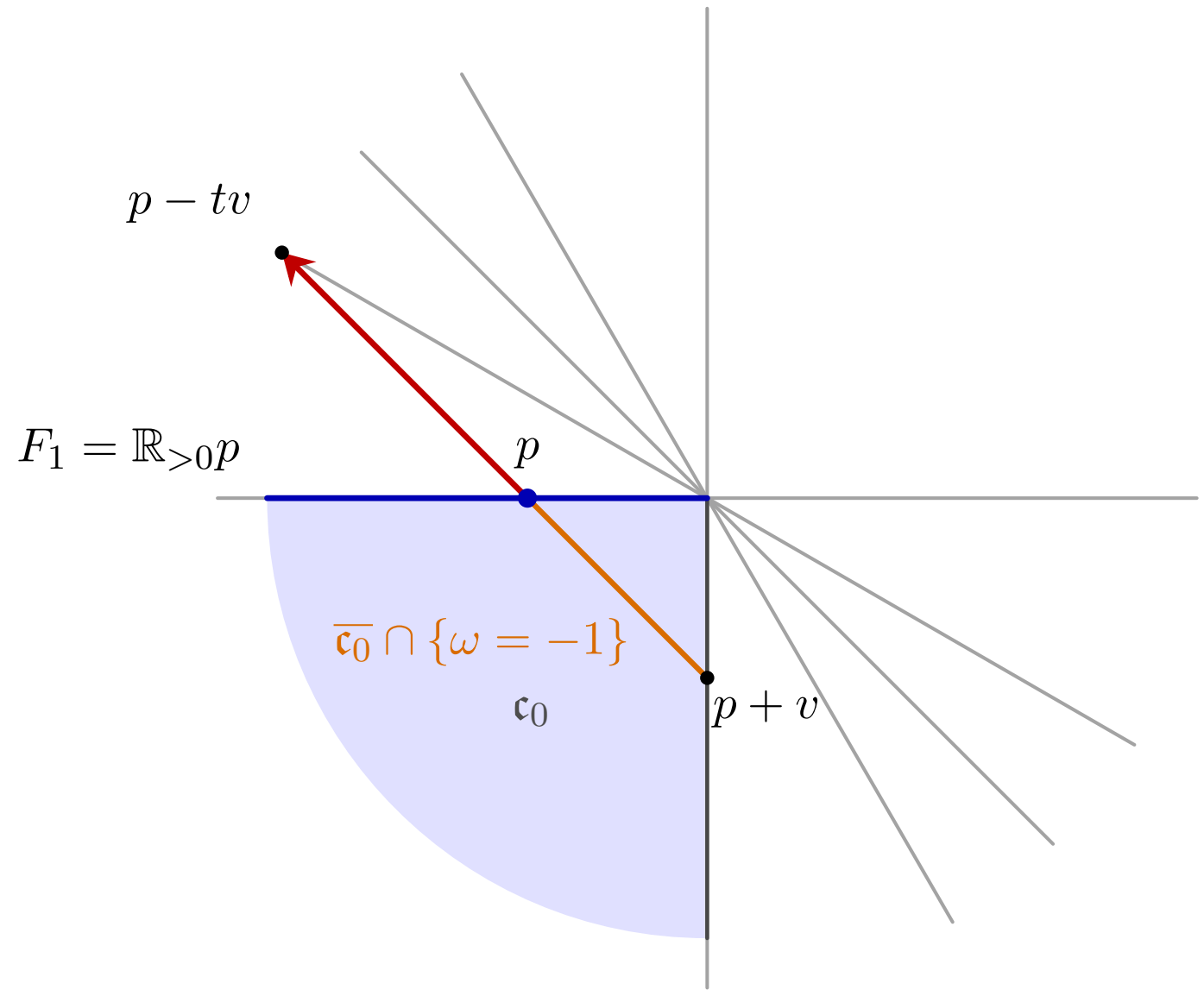
}
\FloatBarrier

\subsection{Flag selection}

In Figure~\ref{fig:two-dimensional-score}, the first slice is cut out by \(\omega=-1\).
At an optimal face in any dimension, \(d\omega\) measures the change in the objective on the normal quotient.

\begin{lemma}\label{lem:direction-cone-structure}
Let \(\mathfrak c\) have finite optimal value and let \(F=F(\mathfrak c,\omega)\). Then
\(d\omega\) induces a covector on \(N_F\) that is negative on \(\overline{R_F(\mathfrak c)}\setminus\{0\}\).
\end{lemma}

\begin{proof}
The objective is constant on \(F\), so its differential vanishes on \(\lin(F)\).
Fix \(x\in F\). Let \(\xi\in\overline{R_F(\mathfrak c)}\setminus\{0\}\) and choose \(v\in V_A\) with \(\pi_F(v)=\xi\).
For sufficiently small \(t>0\), \(x+tv\in\overline{\mathfrak c}\). Since \(\xi\ne0\), this point lies
outside \(\aff(F)\), so \(\omega(x+tv)<\omega(x)\). Thus \(d\omega(v)<0\).
\end{proof}

For each optimal face \(F\), we also write \(d\omega\) for the induced covector
\(N_F\to\R,\ \pi_F(v)\mapsto d\omega(v)\).
Put \(r=\codim F\) and choose an ordered list of covectors on \(N_F\)
\[
\Theta_F=(\theta_{F,r-1},\ldots,\theta_{F,1})
\]
so that \((d\omega,\theta_{F,r-1},\ldots,\theta_{F,1})\) is an ordered basis of \(N_F^*\).
We call this basis the frame at \(F\).
Use the same frame for all chambers with optimal face \(F\).

Start the flag of \(\mathfrak c\) at its optimal face \(F\) with this frame.
Suppose the current face is \(G\subseteq\overline{\mathfrak c}\) of codimension \(i>0\), and the
first covector \(\eta\in N_G^*\) is negative on \(\overline{R_G(\mathfrak c)}\setminus\{0\}\).
At the first step this holds with \(G=F\) and \(\eta=d\omega\) by Lemma~\ref{lem:direction-cone-structure}.
Define the normalized slice by
\[
P_\eta(\mathfrak c,G)
:=
\{\xi\in\overline{R_G(\mathfrak c)}:\eta(\xi)=-1\}.
\]

For each \(H\in\operatorname{Cov}_G(\mathfrak c)\), choose \(\xi_H\in R_G(H)\) and set
\[
p_H:=\frac{\xi_H}{-\eta(\xi_H)}.
\]
The point \(p_H\) is the unique point of the ray \(\overline{R_G(H)}\) on the slice.
By Lemma~\ref{lem:covering-faces-generate-local-directions}, every point in the slice is a nonnegative
combination of the \(p_H\). Applying \(\eta\) shows that the coefficients sum to one, so the slice
is their convex hull. Each \(p_H\) lies on a distinct extreme ray, so these points are precisely its vertices.

When \(i\geq2\), write the current list on \(N_G\) as \((\eta,\varphi,\psi_{i-2},\ldots,\psi_1)\).
When \(i=2\), the list ends with \(\varphi\). Maximize \(\varphi\) on \(P=P_\eta(\mathfrak c,G)\).
If the maximizing vertex \(p\) is unique, take the covering face \(H\) represented by \(p\)
as the next face of the flag. We will show later that we can make the maximizing vertex unique generically.

Let \(\pi_{GH}\colon N_G\to N_H\) be the quotient map, given by
\(\pi_{GH}(\pi_G(v))=\pi_H(v)\) for \(v\in V_A\). Its kernel is \(\lin(H)/\lin(G)=\R p\).
For \(q\in P\), the displacement \(q-p\) lies in \(\ker\eta\) and represents \(\pi_{GH}(q)\).
For each \(\chi\) in \((\varphi,\psi_{i-2},\ldots,\psi_1)\), we therefore want
\(\chi^+(\pi_{GH}(q))=\chi(q-p)=\chi(q)-\chi(p)\).
For arbitrary \(x\in N_G\), the unique representative of \(\pi_{GH}(x)\) in \(\ker\eta\) is
\(x+\eta(x)p\), since \(\eta(p)=-1\). Evaluating \(\chi\) on this representative defines \(\chi^+\in N_H^*\) by
\begin{equation}\label{eq:frame-covector-update}
\chi^+(\pi_{GH}(x)):=\chi\bigl(x+\eta(x)p\bigr)=\chi(x)+\chi(p)\eta(x).
\end{equation}
Take \(\eta^+:=\varphi^+\) as the next normalizing covector. The data at \(H\) are now
\((\eta^+,\psi_{i-2}^+,\ldots,\psi_1^+)\). If \(i>2\), maximize \(\psi_{i-2}^+\) on
\(P_{\eta^+}(\mathfrak c,H)\) to choose the next face.

\begin{figure}[!htbp]
\centering
\IfFileExists{Graph/figures/frame_vertex_figure.tikz}{%
  \input{Graph/figures/frame_vertex_figure.tikz}%
}{%
  \input{figures/frame_vertex_figure.tikz}%
}
\caption{The edge from \(p\) to \(q\) on the left gives the ray through \(q-p\) on the right.
On that ray \(\varphi(v)=-1\) at \(v=(q-p)/(\varphi(p)-\varphi(q))\), the left endpoint of the orange segment.
The next covector chooses between the two endpoints.}
\label{fig:frame-vertex-figure}
\end{figure}
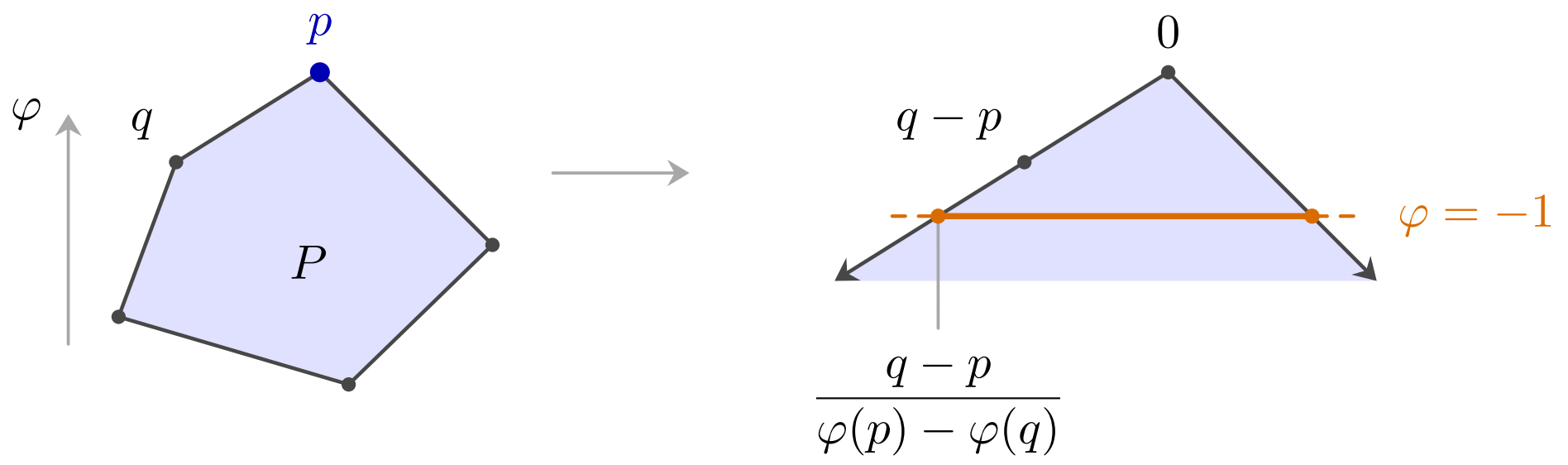

To see the next slice, translate \(p\) to the origin as in Figure~\ref{fig:frame-vertex-figure}.
Extend the edges at \(p\) to rays. On each ray mark the displacement \(v\) for which
\(\varphi(p+v)=\varphi(p)-1\). In the figure there are two such rays. The quotient map \(\pi_{GH}\) sends
the segment joining the marked points onto the next slice \(P_{\eta^+}(\mathfrak c,H)\).

\FloatBarrier

To repeat the step at \(H\), we need \(\eta^+\) to be negative on the nonzero points
in \(\overline{R_H(\mathfrak c)}\).

\begin{lemma}\label{lem:normalized-vertex-figure}
The covector \(\eta^+\) is negative on \(\overline{R_H(\mathfrak c)}\setminus\{0\}\).
The vertices of the next slice \(P_{\eta^+}(\mathfrak c,H)\) are precisely the points
\[
\frac{\pi_{GH}(q-p)}{\varphi(p)-\varphi(q)}
\]
as \(q\) ranges over the vertices of \(P\) joined to \(p\) by an edge.
\end{lemma}

\begin{proof}
Identify \(\ker\eta\) with \(N_H\) by the restriction of \(\pi_{GH}\).

For \(v\in\ker\eta\), we have \(v\in\operatorname{cone}(P-p)\) exactly when \(p+\varepsilon v\in P\) for all
sufficiently small \(\varepsilon>0\). Since \(v\in\ker\eta\), the point \(p+\varepsilon v\) still
lies on \(\eta=-1\). Thus \(p+\varepsilon v\in P\) exactly when it satisfies the inequalities
defining \(\overline{R_G(\mathfrak c)}\). Those that are strict at \(p\) stay strict for small \(\varepsilon\).
Since \(p\in R_G(H)\), the inequalities that hold with equality at \(p\) are precisely those
from hyperplanes containing \(H\). For each such inequality \(a\geq0\), the form \(a\) is linear
on \(N_G\) and \(a(p)=0\), so \(a(p+\varepsilon v)=\varepsilon a(v)\).
These forms descend to \(N_H\), where their inequalities define \(\overline{R_H(\mathfrak c)}\).
Thus \(v\in\operatorname{cone}(P-p)\) exactly when \(\pi_{GH}(v)\in\overline{R_H(\mathfrak c)}\).

The extreme rays of \(\operatorname{cone}(P-p)\) are the rays through \(q-p\) along the edges
of \(P\) incident to \(p\). Their images under \(\pi_{GH}\) are therefore the extreme rays of
\(\overline{R_H(\mathfrak c)}\). On each generator \(\pi_{GH}(q-p)\), the value of \(\eta^+\) is
\(\varphi(q)-\varphi(p)<0\), since \(p\) uniquely maximizes \(\varphi\). These rays generate
the cone, so \(\eta^+\) is negative away from zero. For \(u=\pi_{GH}(q-p)\), we have
\[
(\R_{\geq0}u)\cap P_{\eta^+}(\mathfrak c,H)
=\left\{\frac{u}{-\eta^+(u)}\right\}
=\left\{\frac{\pi_{GH}(q-p)}{\varphi(p)-\varphi(q)}\right\}.
\]
These points are precisely the vertices of the slice, just as for \(P\).
\end{proof}

\begin{definition}\label{def:chamber-to-flag}
Let \(\lambda(\mathfrak c,\omega)<\infty\), and put \(F=F(\mathfrak c,\omega)\) and
\(r=\codim F\). At codimension \(i\), write
\[
\bigl(G_i,\ \eta_i,\ \theta_{F,i-1}^{[i]},\ldots,\theta_{F,1}^{[i]}\bigr)
\]
for the current face \(G_i\) and an ordered list of covectors on \(N_{G_i}\).
The first covector \(\eta_i\) is strictly negative on \(\overline{R_{G_i}(\mathfrak c)}\setminus\{0\}\).
Start at \(G_r:=F\) with \(\eta_r:=d\omega\).
For \(1\leq j<r\), set \(\theta_{F,j}^{[r]}:=\theta_{F,j}\), the corresponding entry of \(\Theta_F\).
At each stage let \(P_i:=P_{\eta_i}(\mathfrak c,G_i)\) be the slice.

For \(i\geq2\), set \(\varphi_i:=\theta_{F,i-1}^{[i]}\). If its values on the vertices of \(P_i\)
are pairwise distinct, apply the step above with \(\varphi=\varphi_i\) and
\(\psi_j=\theta_{F,j}^{[i]}\) for \(1\leq j\leq i-2\). Denote the selected vertex by \(p_i\).
Write
\[
\bigl(G_{i-1},\ \eta_{i-1},\ \theta_{F,i-2}^{[i-1]},\ldots,\theta_{F,1}^{[i-1]}\bigr)
\]
for the resulting face \(H\) and the covectors \((\eta^+,\psi_{i-2}^+,\ldots,\psi_1^+)\)
of~\eqref{eq:frame-covector-update}. Repeat at \(i-1\).

If the recursion reaches \(i=1\), only \(\eta_1\) remains and the only face covering \(G_1\) in
\(\overline{\mathfrak c}\) is \(\mathfrak c\) itself. Set \(G_0=\mathfrak c\) and define
\[
\mathcal F_{\mathrm F}(\mathfrak c):=(G_r,G_{r-1},\ldots,G_0=\mathfrak c).
\]
If \(\lambda(\mathfrak c,\omega)=\infty\), choose the codimension-zero flag
\(\mathcal F_{\mathrm F}(\mathfrak c):=(\mathfrak c)\).

Call \(\Theta_F\) admissible at \(F\) if, in every chamber with optimal face \(F\),
the covector \(\varphi_i\) takes distinct values on the vertices of \(P_i\) at every step with \(i \geq 2\), so the recursion is well defined.
A family \(\Theta=(\Theta_F)_F\) is admissible if \(\Theta_F\) is admissible at each optimal face \(F\).
For an admissible family, let \(G_{\mathrm F}:=G_{\mathcal F_{\mathrm F}}\) be the graph of
the flag selection strategy \(\mathcal F_{\mathrm F}\) as in Definition~\ref{def:graph-from-flag-strategy}.
\end{definition}

\begin{lemma}\label{lem:generic-regular-face-frames}\label{lem:regular-frame-system-exists}
For each optimal face \(F\) of codimension \(r\), the set of admissible choices of \(\Theta_F\)
contains a nonempty Zariski open subset of \((N_F^*)^{r-1}\).
In particular, an admissible family \(\Theta=(\Theta_F)_F\) exists.
\end{lemma}

\begin{proof}
Fix \(F\). Along any sequence of faces reached before a tie, each slice vertex has the form
\(-\xi_H/\eta(\xi_H)\) for a fixed ray generator \(\xi_H\). Starting with \(\eta=d\omega\),
the update~\eqref{eq:frame-covector-update} therefore makes the vertices and their covector values
rational functions of the entries of \(\Theta_F\).

For distinct current vertices \(u,v\), the difference \(\varphi_i(u)-\varphi_i(v)\) is a nonzero rational function.
Indeed, fix the earlier choices and take a covector on the current quotient separating \(u\) and \(v\).
Use its pullback to \(N_F\) as the next frame entry. At each quotient step the pullback vanishes
on the kernel, so~\eqref{eq:frame-covector-update} recovers the chosen covector on the current quotient.

Only finitely many chambers and sequences of faces occur at \(F\), by Lemma~\ref{lem:star-arrangement-correspondence}.
Avoiding the zeros of the numerators and denominators of all these differences, and requiring
\((d\omega,\Theta_F)\) to be a basis, gives a nonempty Zariski open set of admissible choices.
\end{proof}

\subsection{Strict increase of the potential}\label{sec:strategy-proof}

Fix an admissible family \(\Theta\). For a chamber with finite optimal value, write
\[
S_i^{\mathrm F}(\mathfrak c):=\varphi_i(p_i)=\max_{q\in P_i}\varphi_i(q)
\qquad(2\leq i\leq\codim F(\mathfrak c,\omega)).
\]

\begin{definition}\label{def:frame-potential}
If \(\lambda(\mathfrak c,\omega)<\infty\), put \(r=\codim F(\mathfrak c,\omega)\) and define the potential by
\[
P_{\mathrm F}^*(\mathfrak c)
:=
\bigl(\lambda(\mathfrak c,\omega),\,
      \dim F(\mathfrak c,\omega),\,
      S_r^{\mathrm F}(\mathfrak c),\ldots,S_2^{\mathrm F}(\mathfrak c)\bigr).
\]
The last face of the flag is forced to be \(\mathfrak c\), so the sequence stops at codimension two.

Compare these tuples lexicographically.
If their first two coordinates agree, their optimal faces have the same codimension, so the
remaining tuples have the same length.

If \(\lambda(\mathfrak c,\omega)=\infty\), set \(P_{\mathrm F}^*(\mathfrak c)=+\infty\),
ordered above every chamber with finite optimal value.
\end{definition}

\begin{proposition}\label{prop:first-changed-face-better-target-flag}\prooflink{prop:first-changed-face-better-target-flag}
For every edge \(\mathfrak c\to\mathfrak b\) in \(G_{\mathrm F}\),
\[
P_{\mathrm F}^*(\mathfrak b)>_{\mathrm{lex}}P_{\mathrm F}^*(\mathfrak c).
\]
\end{proposition}

\begin{proof}[\prooftarget{prop:first-changed-face-better-target-flag}{Proof of Proposition~\ref*{prop:first-changed-face-better-target-flag}}]
The source has finite optimal value, since chambers with infinite optimal value have no outgoing edges.
By Lemma~\ref{lem:flag-strategy-weak-increase},
\(\lambda(\mathfrak b,\omega)\geq\lambda(\mathfrak c,\omega)\).
If the target has infinite optimal value, its potential is \(+\infty\), so the conclusion holds.
If both optimal values are finite and the inequality is strict, the first coordinate of the potential increases.
Otherwise the same lemma gives
\[
\overline{F(\mathfrak c,\omega)}\subseteq\overline{F(\mathfrak b,\omega)}.
\]
If the optimal faces have different dimensions, the second coordinate increases.
We may therefore assume they are the same face \(F\), of codimension \(r\).

Write \(\mathcal F_{\mathrm F}(\mathfrak c)=(F_r,\ldots,F_0=\mathfrak c)\) and
\(\mathcal F_{\mathrm F}(\mathfrak b)=(G_r,\ldots,G_0=\mathfrak b)\). By
Definition~\ref{def:graph-from-flag-strategy}, there is \(\alpha\in\{\pm1\}^r\setminus\{1^r\}\) such that
\(\mathfrak b=\mathfrak b_{\mathcal F_{\mathrm F}(\mathfrak c)}(\alpha)\).
Besides the two selected flags, consider the flag \((F_r(\alpha),\ldots,F_0(\alpha))\) given by
Theorem~\ref{thm:opposite-flag}. Put \(i^*:=\max\{i:\alpha_i=-1\}\).
Lemma~\ref{lem:first-flag-change} gives
\[
F_i(\alpha)=F_i\quad(i^*\leq i\leq r),\qquad
R_{F_{i^*}}(F_{i^*-1}(\alpha))=-R_{F_{i^*}}(F_{i^*-1}).
\]

We compare the two selected flags starting from their common base \(F_r=G_r=F\), so both
chambers use \(\eta_r=d\omega\) and the covectors in \(\Theta_F\) in the same order.
Whenever they choose the same covering face from the same data, its ray meets \(\{\eta_i=-1\}\)
at the same point \(p_i\). The updates therefore give the same covectors on the common quotient
at the next face.

Suppose first that they make different choices at some level \(j>i^*\), and take the first
such level in the descending recursion. Thus
\[
F_a=G_a\quad(j\leq a\leq r),\qquad F_{j-1}\neq G_{j-1}.
\]
Since \(j-1\geq i^*\), Lemma~\ref{lem:first-flag-change} gives
\(F_{j-1}=F_{j-1}(\alpha)\subseteq\overline{\mathfrak b}\).
It is therefore an available covering face at \(G_j=F_j\) in the target's selection.
Its ray meets \(\{\eta_j=-1\}\) at the same point in both slices, so \(\varphi_j\) takes the
value \(S_j^{\mathrm F}(\mathfrak c)\) there. The target selects a different vertex, its unique maximizer, so
\[
S_j^{\mathrm F}(\mathfrak b)>S_j^{\mathrm F}(\mathfrak c).
\]
All earlier maxima \(S_a^{\mathrm F}\), \(a>j\), agree. Hence the target potential is strictly larger.

It remains to consider the case in which the selected flags agree through \(F_{i^*}\).
They then have the same covector \(\eta_{i^*}\) there. The cone \(\overline{R_{F_{i^*}}(\mathfrak c)}\)
contains the ray into \(F_{i^*-1}\), and \(\overline{R_{F_{i^*}}(\mathfrak b)}\) contains the opposite
ray into \(F_{i^*-1}(\alpha)\). The covector \(\eta_{i^*}\) is negative away from
zero on both cones, which is impossible for opposite rays.
\end{proof}

\begin{theorem}\label{theorem:acyclic}\label{theorem:regular-frame-win-conditions}\label{theorem:regular-quotient-lex-win-conditions}\prooflink{theorem:regular-frame-win-conditions}
Let \(\mathcal P\) be a locally finite affine hyperplane arrangement in \(A\), and let
\(\omega\colon A\to\R\) be a nonconstant affine objective.
At each optimal face \(F\), there exists an admissible list of covectors \(\Theta_F\).
For every such family \(\Theta=(\Theta_F)_F\), the graph \(G_{\mathrm F}\) is acyclic and the strategy
\(\mathcal F_{\mathrm F}\) guarantees both progress along every infinite walk and uniform progress. If \(\mathcal P\) has finitely many
normal spaces, then \(\mathcal F_{\mathrm F}\) also guarantees both escape along every infinite walk and uniform escape.
\end{theorem}

\begin{proof}[\prooftarget{theorem:regular-frame-win-conditions}{Proof of Theorem~\ref*{theorem:regular-frame-win-conditions}}]
Lemma~\ref{lem:regular-frame-system-exists} gives an admissible list at each optimal face. For any such choice,
Proposition~\ref{prop:first-changed-face-better-target-flag} makes the potential strictly increase along
every directed edge, which excludes directed cycles. The progress and escape assertions now follow from
Corollary~\ref{cor:acyclic-flag-strategy-progressive} and Theorem~\ref{thm:acyclic-flag-strategy-escaping}.
\end{proof}

\section*{Use of generative AI}

This paper combines results developed by the author with results developed by ChatGPT.
In fall 2025, the author developed the construction of a graph traversal from flags, the Euclidean strategy,
the proof of its acyclicity, the proof that acyclicity implies the progress conditions, and the
proof that finitely many normal spaces together with acyclicity imply the escape conditions.
ChatGPT was used as a discussion partner in developing these results. Most remaining results
were developed primarily by ChatGPT in summer 2026, with guidance and suggestions from the
author, followed by substantial human review, editing, and reformulation.

\end{document}

%% file: figures/nonpolyhedral_chamber.tikz
\begin{tikzpicture}[x=0.98cm,y=1.05cm,line cap=round,line join=round,
  boundary/.style={draw=blue!68!black,line width=1.35pt},
  wall/.style={draw=orange!80!black,line width=0.55pt,dash pattern=on 1.1pt off 2.7pt},
  point/.style={circle,fill=black,inner sep=1.0pt},
  smalllabel/.style={font=\scriptsize,inner sep=2pt}]
  \foreach \n in {0,...,7}
    \coordinate (q\n) at (\n,{sqrt(\n)});

  \begin{scope}
    \clip (-0.48,-0.60) rectangle (6.65,3.85);
    \fill[blue!10] (-0.48,-0.65) -- (7,-0.65) -- (q7)
      \foreach \n in {6,...,0} { -- (q\n) } -- (-0.48,-0.48) -- cycle;
    \foreach \n in {0,...,6} {
      \pgfmathsetmacro{\slope}{1/(sqrt(\n+1)+sqrt(\n))}
      \pgfmathsetmacro{\intercept}{sqrt(\n)-\n*\slope}
      \draw[wall] (-0.48,{\intercept-0.48*\slope}) -- (6.65,{\intercept+6.65*\slope});
    }
    \draw[boundary] (-0.48,-0.48) -- (q0) \foreach \n in {1,...,7} { -- (q\n) };
    \foreach \n in {0,...,6} \node[point] at (q\n) {};
  \end{scope}
  \node[smalllabel,below right] at (q0) {$q_0$};
  \node[smalllabel,below right] at (q1) {$q_1$};
  \node[smalllabel,below right] at (q4) {$q_4$};
\end{tikzpicture}

%% file: figures/unattained_optimum.tikz
\begin{tikzpicture}[x=0.72cm,y=3.0cm,line cap=round,line join=round,
  boundary/.style={draw=blue!68!black,line width=1.35pt},
  wall/.style={draw=orange!80!black,line width=0.55pt,dash pattern=on 1.1pt off 2.7pt},
  guide/.style={draw=black!34,line width=0.7pt,dash pattern=on 1.1pt off 2.7pt},
  point/.style={circle,fill=black,inner sep=1.0pt}]
  \foreach \n in {0,...,10}
    \coordinate (p\n) at (\n,{-1/(2^\n)});
  \coordinate (left) at (-1,-1.5);

  \begin{scope}
    \clip (-1,-1.54) rectangle (9.25,0.11);
    \fill[blue!10] (-1,-1.54) -- (10,-1.54) -- (p10)
      \foreach \n in {9,...,0} { -- (p\n) } -- (left) -- cycle;
    \foreach \n in {0,...,9} {
      \pgfmathsetmacro{\slope}{1/(2^(\n+1))}
      \pgfmathsetmacro{\intercept}{-(\n+2)*\slope}
      \draw[wall] (-1,{\intercept-\slope}) -- (9.25,{\intercept+9.25*\slope});
    }
    \draw[guide] (-1,0) -- (9.25,0);
    \draw[boundary] (left) \foreach \n in {0,...,10} { -- (p\n) };
    \foreach \n in {0,...,9} \node[point] at (p\n) {};
  \end{scope}
  \node[right,font=\scriptsize] at (9.25,0) {$y=0$};
  \node[blue!65!black,font=\scriptsize] at (6.7,-0.62) {$C$};
\end{tikzpicture}

%% file: figures/attainment_fixed_segment.tikz
\begin{tikzpicture}[line cap=round, line join=round, >={Stealth[length=7pt,width=8pt]},
  boundary/.style={draw=blue!68!black,line width=1.35pt},
  support/.style={draw=orange!85!black,line width=0.85pt,dash pattern=on 1.0pt off 2.6pt},
  guide/.style={draw=black!40,line width=0.85pt,dash pattern=on 1.0pt off 2.6pt},
  point/.style={circle,fill=black,inner sep=1.15pt},
  hit/.style={circle,fill=orange!85!black,inner sep=1.35pt}]
  \coordinate (P) at (0,-0.90);
  \coordinate (S) at (0,2.05);
  \coordinate (Intersection) at (0,0.63);
  \coordinate (Z) at (3.20,1.10);
  \coordinate (Xn) at ($(S)!1.52!(Z)$);
  \coordinate (Bend) at ($(Intersection)!0.55!(Z)$);
  \coordinate (Tail) at ($(Intersection)!1.8125!(Z)$);
  \coordinate (Left) at (-0.55,-1.00);

  \fill[blue!9] (-0.55,-1.30) -- (5.80,-1.30) -- (Tail) -- (Bend) -- (Left) -- cycle;
  \node[blue!65!black] at (5.12,-0.72) {$C$};
  \draw[boundary] (Left) -- (Bend) -- (Tail);
  \draw[support] ($(Intersection)!-0.15!(Z)$) -- (Tail);
  \node[orange!80!black,above] at ($(Intersection)!0.32!(Z)$) {$H_n$};

  \draw[guide] (-0.35,2.05) -- (5.80,2.05) node[right,black] {$\omega=0$};
  \draw[black!65,dashed] (Xn) -- (S);
  \draw[orange!85!black,line width=2.0pt] (P) -- (S);
  \draw[->,black!58,line width=0.7pt] (-0.74,-0.95) -- (-0.74,-0.32) node[left] {$e$};

  \node[point,label=right:$p$] at (P) {};
  \node[hit] at (Intersection) {};
  \node[hit,label=above right:$z_n$] at (Z) {};
  \node[point,label=below right:$x_n$] at (Xn) {};
  \node[hit,label=above left:$s$] at (S) {};
\end{tikzpicture}

%% file: figures/flag_choice_perturbations.tikz
\begin{tikzpicture}[
  scale=1.0,
  line cap=round,
  line join=round,
  >={Stealth[length=6pt,width=7pt]},
  arrangement/.style={draw=gray!72,line width=0.75pt},
  start/.style={fill=blue!30,opacity=0.82},
  flagface/.style={draw=blue!70!black,line width=1.65pt},
  perturb/.style={draw=orange!85!black,line width=0.55pt},
  sourcepoint/.style={fill=black,draw=white,line width=0.35pt},
  targetpoint/.style={fill=orange!85!black,draw=white,line width=0.35pt},
  paneltitle/.style={font=\small},
  pointlabel/.style={font=\scriptsize,inner sep=1.2pt}
]

\newcommand{\flagchoicepanel}[4]{%
  \begin{scope}[shift={(#1,#2)}]
    \def\R{2.05}
    \def\L{2.34}
    \path[start] (0,0) -- (180:\R) arc[start angle=180,end angle=270,radius=\R] -- cycle;
    \draw[arrangement] (-\L,0) -- (\L,0);
    \draw[arrangement] (0,\L) -- (0,-\L);
    \draw[arrangement] (-2.03,1.17) -- (2.03,-1.17);
    \draw[arrangement] (-1.65,1.65) -- (1.65,-1.65);
    \draw[arrangement] (-1.17,2.03) -- (1.17,-2.03);
    #3
    \filldraw[sourcepoint] #4
  \end{scope}
}

\flagchoicepanel{-2.85}{0}{%
  \draw[flagface] (0,0) -- (-2.05,0);
  \node[paneltitle] at (0,2.55) {$F_1=\{(x_1,0):x_1<0\}$};
  \draw[perturb] (0,0) -- (-1.13,-0.27);
  \draw[perturb] (0,0) -- (-1.13,0.27);
  \draw[perturb] (0,0) -- (1.13,-0.27);
  \draw[perturb] (0,0) -- (1.13,0.27);
  \filldraw[targetpoint] (-1.13,0.27) circle (0.055) node[above left=0pt,pointlabel] {$(-,+)$};
  \filldraw[targetpoint] (1.13,-0.27) circle (0.055) node[below right=0pt,pointlabel] {$(+,-)$};
  \filldraw[targetpoint] (1.13,0.27) circle (0.055) node[above right=0pt,pointlabel] {$(-,-)$};
}{%
  (-1.13,-0.27) circle (0.06) node[below left=0pt,pointlabel] {$(+,+)$};
}

\flagchoicepanel{2.85}{0}{%
  \draw[flagface] (0,0) -- (0,-2.05);
  \node[paneltitle] at (0,2.55) {$F_1=\{(0,x_2):x_2<0\}$};
  \draw[perturb] (0,0) -- (-0.27,-1.13);
  \draw[perturb] (0,0) -- (0.27,-1.13);
  \draw[perturb] (0,0) -- (-0.27,1.13);
  \draw[perturb] (0,0) -- (0.27,1.13);
  \filldraw[targetpoint] (0.27,-1.13) circle (0.055) node[below right=0pt,pointlabel] {$(-,+)$};
  \filldraw[targetpoint] (-0.27,1.13) circle (0.055) node[above left=0pt,pointlabel] {$(+,-)$};
  \filldraw[targetpoint] (0.27,1.13) circle (0.055) node[above right=0pt,pointlabel] {$(-,-)$};
}{%
  (-0.27,-1.13) circle (0.06) node[below left=0pt,pointlabel] {$(+,+)$};
}

\end{tikzpicture}

%% file: figures/traversal_example.tikz
\begin{tikzpicture}[scale=0.94, line cap=round, line join=round, >={Stealth[length=7pt,width=8pt]}]
  \path[fill=orange!28] (0,0) -- (0:4.050) arc[start angle=0, end angle=90, radius=4.050] -- cycle;
  \path[fill=orange!28] (0,0) -- (90:4.050) arc[start angle=90, end angle=120, radius=4.050] -- cycle;
  \path[fill=orange!28] (0,0) -- (120:4.050) arc[start angle=120, end angle=135, radius=4.050] -- cycle;
  \path[fill=blue!12] (0,0) -- (135:4.050) arc[start angle=135, end angle=150, radius=4.050] -- cycle;
  \path[fill=blue!12] (0,0) -- (150:4.050) arc[start angle=150, end angle=180, radius=4.050] -- cycle;
  \path[fill=blue!12] (0,0) -- (180:4.050) arc[start angle=180, end angle=270, radius=4.050] -- cycle;
  \path[fill=blue!12] (0,0) -- (270:4.050) arc[start angle=270, end angle=300, radius=4.050] -- cycle;
  \path[fill=blue!12] (0,0) -- (300:4.050) arc[start angle=300, end angle=315, radius=4.050] -- cycle;
  \path[fill=orange!28] (0,0) -- (315:4.050) arc[start angle=315, end angle=330, radius=4.050] -- cycle;
  \path[fill=orange!28] (0,0) -- (330:4.050) arc[start angle=330, end angle=360, radius=4.050] -- cycle;
  \path[fill=blue!32, opacity=0.75] (0,0) -- (180:4.050) arc[start angle=180, end angle=270, radius=4.050] -- cycle;
  \draw[draw=gray!78, line width=0.9pt] (-4.600, 0.000) -- (4.600, 0.000);
  \draw[draw=gray!78, line width=0.9pt] (-3.984, 2.300) -- (3.984, -2.300);
  \draw[draw=gray!78, line width=0.9pt] (-3.253, 3.253) -- (3.253, -3.253);
  \draw[draw=gray!78, line width=0.9pt] (-2.300, 3.984) -- (2.300, -3.984);
  \draw[draw=gray!78, line width=0.9pt] (0.000, 4.600) -- (0.000, -4.600);
  \draw[draw=orange!90!black, line width=1.0pt, dotted, postaction={decorate}, decoration={markings, mark=at position 0.95 with {\arrow{Stealth[length=7pt,width=8pt]}}}] (225:4.600) -- (45:4.600);
  \draw[->, draw=black!75, line width=0.9pt] (-1.770, -1.770) -- (-2.418, 0.648);
  \draw[->, draw=black!75, line width=0.9pt] (-1.770, -1.770) -- (1.860, 1.675);
  \draw[->, draw=black!75, line width=0.9pt] (-1.770, -1.770) -- (2.418, -0.648);
  \draw[->, draw=black!75, line width=0.9pt] (-2.418, 0.648) -- (-1.986, 1.524);
  \draw[->, draw=black!75, line width=0.9pt] (-2.418, 0.648) -- (1.986, -1.524);
  \draw[->, draw=black!75, line width=0.9pt] (-2.418, 0.648) -- (2.418, -0.648);
  \draw[->, draw=black!75, line width=0.9pt] (-1.986, 1.524) -- (-1.524, 1.986);
  \filldraw[fill=black, draw=white, line width=0.4pt] (-1.770, -1.770) circle (0.075);
  \filldraw[fill=black, draw=white, line width=0.4pt] (-2.418, 0.648) circle (0.075);
  \filldraw[fill=black, draw=white, line width=0.4pt] (1.860, 1.675) circle (0.075);
  \filldraw[fill=black, draw=white, line width=0.4pt] (2.418, -0.648) circle (0.075);
  \filldraw[fill=black, draw=white, line width=0.4pt] (-1.986, 1.524) circle (0.075);
  \filldraw[fill=black, draw=white, line width=0.4pt] (1.986, -1.524) circle (0.075);
  \filldraw[fill=black, draw=white, line width=0.4pt] (-1.524, 1.986) circle (0.075);
\end{tikzpicture}

%% file: figures/unit_circle_selection.tikz
\begin{tikzpicture}[
  x=2cm,
  y=2cm,
  line cap=round,
  line join=round,
  >={Stealth[length=5pt,width=6pt]},
  axis/.style={draw=gray!65,line width=0.65pt},
  unitcircle/.style={draw=gray!60,densely dotted,line width=0.8pt},
  candidatearc/.style={draw=blue!60!black,line width=1.1pt},
  objectivelevel/.style={draw=orange!85!black,dashed,line width=0.85pt},
  generator/.style={draw=black!72,line width=1pt},
  chosen/.style={draw=red!70!black,line width=1.35pt}
]
  \coordinate (origin) at (0,0);
  \coordinate (leftGenerator) at (-1,0);
  \coordinate (lowerGenerator) at (0,-1);

  \path[fill=blue!12] (origin) -- (leftGenerator)
    arc[start angle=180,end angle=270,radius=2cm] -- cycle;
  \foreach \angle in {120,135,150}
    \draw[axis] (\angle:1.22) -- ({\angle+180}:1.22);
  \draw[->,axis] (-1.23,0) -- (1.18,0) node[right=2pt] {$x_1$};
  \draw[->,axis] (0,-1.23) -- (0,1.18) node[above=2pt] {$x_2$};
  \draw[unitcircle] (origin) circle[radius=2cm];
  \draw[candidatearc] (leftGenerator) arc[start angle=180,end angle=270,radius=2cm];
  \node[above right=2pt,text=black!70] at (49:1) {$x_1^2+x_2^2=1$};
  \node[text=black!70] at (-0.28,-0.28) {$\mathfrak c_0$};

  \draw[objectivelevel] (-1.25,0.25) -- (0.18,-1.18);
  \node[above left=1pt,text=orange!85!black] at (-1.25,0.25) {$\ell=-1$};
  \draw[->,generator] (origin) -- (lowerGenerator);
  \draw[->,chosen] (origin) -- (leftGenerator);
  \fill[black!72] (origin) circle[radius=1.1pt];
  \fill[black!72] (lowerGenerator) circle[radius=1.4pt];
  \fill[red!70!black] (leftGenerator) circle[radius=1.65pt];
  \node[below right=3pt] at (lowerGenerator) {$(0,-1)$};
  \node[below left=3pt,text=red!70!black,align=right] at (leftGenerator)
    {$(-1,0)$\\[-1pt]\scriptsize selected};
\end{tikzpicture}

%% file: figures/two_intermediate_faces.tikz
\begin{tikzpicture}[
  line cap=round,
  line join=round,
  >={Stealth[length=6pt,width=7pt]},
  chamber/.style={fill=blue!12},
  ray/.style={draw=black!72,line width=0.9pt},
  chosenray/.style={draw=blue!70!black,line width=1.25pt},
  generator/.style={draw=orange!85!black,line width=1.4pt},
  normal/.style={draw=red!75!black,line width=1.25pt},
  component/.style={draw=orange!85!black,dashed,line width=0.8pt},
  unitcircle/.style={draw=gray!60,densely dotted,line width=0.7pt},
  point/.style={circle,fill=black,inner sep=1.1pt}
]
  \foreach \shift/\angle/\side in {0/55/right,7.2/125/left}{
    \begin{scope}[xshift=\shift cm]
      \coordinate (O) at (0,0);
      \coordinate (U) at (1.75,0);
      \coordinate (V) at (1.75,1.75);
      \coordinate (W) at (\angle:1.75);
      \coordinate (A) at ({1.75*cos(\angle)},0);

      \path[chamber] (O) -- (3.2,0)
        arc[start angle=0,end angle=\angle,radius=3.2] -- cycle;
      \draw[gray!55,line width=0.65pt] (-1.95,0) -- (3.45,0);
      \draw[unitcircle] (1.75,0)
        arc[start angle=0,end angle=180,radius=1.75];
      \draw[ray] (O) -- (\angle:3.35)
        node[above=3pt] {$R_F(H')$};
      \draw[chosenray] (O) -- (3.45,0)
        node[below=5pt] {$R_F(H)$};

      \draw[component] (A) -- (W)
        node[pos=0.65,\side=3pt,orange!85!black] {$b v$};
      \draw[->,chosenray,line width=1.55pt] (O) -- (U)
        node[below=5pt,blue!70!black] {$u$};
      \draw[->,normal] (U) -- (V)
        node[above=4pt,red!75!black] {$v$};
      \draw[->,generator] (O) -- (W);
      \ifnum\angle<90
        \node[above left=2pt,orange!85!black] at (W) {$w$};
      \else
        \node[above right=2pt,orange!85!black] at (W) {$w$};
      \fi

      \draw[gray!70,line width=0.7pt] (1.94,0) -- (1.94,0.19) -- (1.75,0.19);
      \node[point] at (O) {};
      \node[below=17pt] at (O) {$R_F(F)=0$};
      \node[point] at (A) {};
      \node[below=5pt] at (A) {$a u$};
      \node[black!70] at (2.5,0.6) {$R_F(K)$};
    \end{scope}
  }
\end{tikzpicture}

%% file: figures/two_dimensional_score.tikz
\begin{tikzpicture}[scale=1.0, line cap=round, line join=round, >={Stealth[length=7pt,width=8pt]},
  chamber/.style={fill=blue!12},
  arrangement/.style={draw=gray!72,line width=0.8pt},
  ray/.style={draw=black!72,line width=0.9pt},
  chosenray/.style={draw=blue!70!black,line width=1.35pt},
  slice/.style={draw=orange!85!black,line width=1.25pt},
  tangent/.style={draw=blue!70!black,line width=1.55pt},
  stretch/.style={draw=red!75!black,line width=1.35pt},
  point/.style={circle,fill=black,inner sep=1.15pt},
  chosenpoint/.style={circle,fill=blue!70!black,inner sep=1.55pt}]
  \coordinate (O) at (0,0);
  \coordinate (p) at (-1.45,0);
  \coordinate (q) at (0,-1.45);
  \coordinate (r) at (150:{1.45/(cos(30)-sin(30))});
  \coordinate (xext) at (-3.55,0);
  \coordinate (yext) at (0,-3.55);

  \path[chamber] (O)--(180:3.55) arc[start angle=180,end angle=270,radius=3.55] -- cycle;

  \draw[arrangement] (-3.95,0) -- (3.95,0);
  \draw[arrangement] (0,3.95) -- (0,-3.95);
  \draw[arrangement] (150:3.98) -- (-30:3.98);
  \draw[arrangement] (-2.79,2.79) -- (2.79,-2.79);
  \draw[arrangement] (-1.98,3.42) -- (1.98,-3.42);

  \draw[ray] (O)--(yext);
  \draw[chosenray] (O)--(xext) node[above left=2pt] {$F_1=\R_{>0}p$};
  \node[black!70] at (-1.42,-1.72) {$\mathfrak c_0$};

  \draw[slice] (p)--(q);
  \node[orange!85!black] at (-1.82,-1.16) {$\overline{\mathfrak c_0}\cap\{\omega=-1\}$};

  \draw[->,stretch] (p)--(r);

  \node[point] at (q) {};
  \node[anchor=west] at (-0.09,-1.70) {$p+v$};
  \node[point] at (r) {};
  \node[above left=3pt] at (r) {$p-tv$};
  \node[chosenpoint] at (p) {};
  \node[above=3pt] at (p) {$p$};
\end{tikzpicture}

%% file: figures/frame_vertex_figure.tikz
\begin{tikzpicture}[
  scale=1.15,
  line cap=round,
  line join=round,
  >={Stealth[length=5pt,width=6pt]},
  boundary/.style={draw=black!72,line width=0.85pt},
  objective/.style={draw=gray!68,line width=0.7pt},
  level/.style={draw=orange!85!black,dashed,line width=0.8pt},
  newslice/.style={draw=orange!85!black,line width=1.4pt},
  point/.style={circle,fill=black!72,inner sep=1.15pt},
  chosen/.style={circle,fill=blue!70!black,inner sep=1.6pt}
]
  \coordinate (p) at (0,1.25);
  \coordinate (q) at (-1,0.625);
  \coordinate (r) at (1.2,0.05);
  \coordinate (s) at (0.2,-0.92);
  \coordinate (t) at (-1.4,-0.45);

  \path[fill=blue!12, boundary] (p) -- (q) -- (t) -- (s) -- (r) -- cycle;
  \draw[->,objective] (-1.75,-0.64) -- (-1.75,0.96)
    node[left=2pt] {$\varphi$};
  \node at (-0.08,-0.07) {$P$};
  \node[chosen,label={[text=blue!70!black]above:$p$}] at (p) {};
  \node[point,label=above left:$q$] at (q) {};
  \node[point] at (r) {};
  \node[point] at (s) {};
  \node[point] at (t) {};

  \draw[->,objective] (1.62,0.55) -- (2.57,0.55);

  \begin{scope}[shift={(5.9,1.25)}]
    \coordinate (origin) at (0,0);
    \coordinate (qDirection) at (-1,-0.625);
    \coordinate (leftEndpoint) at (-1.6,-1);
    \coordinate (rightEndpoint) at (1,-1);
    \coordinate (leftRayEnd) at (-2.32,-1.45);
    \coordinate (rightRayEnd) at (1.45,-1.45);

    \path[fill=blue!12] (origin) -- (leftRayEnd) -- (rightRayEnd) -- cycle;
    \draw[->,boundary] (origin) -- (leftRayEnd);
    \draw[->,boundary] (origin) -- (rightRayEnd);
    \draw[level] (-1.92,-1) -- (1.38,-1)
      node[right=3pt,text=orange!85!black] {$\varphi=-1$};
    \draw[newslice] (leftEndpoint) -- (rightEndpoint);

    \node[point,label=above:$0$] at (origin) {};
    \node[point,label=above left:$q-p$] at (qDirection) {};
    \draw[objective] (leftEndpoint) -- (-1.6,-1.78);
    \fill[orange!85!black] (leftEndpoint) circle[radius=1.5pt];
    \fill[orange!85!black] (rightEndpoint) circle[radius=1.5pt];
    \node[below=1pt] at (-1.6,-1.78)
      {$\dfrac{q-p}{\varphi(p)-\varphi(q)}$};
  \end{scope}
\end{tikzpicture}

%% file: Graph.bbl
\begin{thebibliography}{99}

\bibitem{AguiarMahajanArrangements}
Marcelo Aguiar and Swapneel Mahajan.
\newblock {\em Topics in Hyperplane Arrangements}, volume 226 of {\em
  Mathematical Surveys and Monographs}.
\newblock American Mathematical Society, Providence, RI, 2017.
\newblock \href {https://doi.org/10.1090/surv/226}
  {\path{doi:10.1090/surv/226}}.
\newblock \href{https://www.math.iitb.ac.in/~swapneel/b.pdf}{Authors' PDF}.

\bibitem{CernyHladikRada2019}
Michal \v{C}ern\'y, Milan Hlad\'ik, and Miroslav Rada.
\newblock Walks on hyperplane arrangements and optimization of piecewise linear functions.
\newblock \href{https://arxiv.org/abs/1912.12750}{arXiv:1912.12750}, 2019.

\bibitem{DeConciniProcesiArrangements}
Corrado~De Concini and Claudio Procesi.
\newblock {\em Topics in Hyperplane Arrangements, Polytopes and Box-Splines}.
\newblock Universitext. Springer, New York, 2010.
\newblock \href {https://doi.org/10.1007/978-0-387-78963-7}
  {\path{doi:10.1007/978-0-387-78963-7}}.

\bibitem{DiestelGraphTheory}
Reinhard Diestel.
\newblock {\em Graph Theory}, volume 173 of {\em Graduate Texts in
  Mathematics}.
\newblock Springer, Berlin, 5th edition, 2017.
\newblock \href {https://doi.org/10.1007/978-3-662-53622-3}
  {\path{doi:10.1007/978-3-662-53622-3}}.

\bibitem{Edelman1984}
Paul~H. Edelman.
\newblock A partial order on the regions of \(\mathbb R^n\) dissected by hyperplanes.
\newblock {\em Transactions of the American Mathematical Society}, 283(2):617--631, 1984.
\newblock \href {https://doi.org/10.1090/S0002-9947-1984-0737888-6}
  {\path{doi:10.1090/S0002-9947-1984-0737888-6}}.

\bibitem{GrunbaumConvexPolytopes}
Branko Gr\"unbaum.
\newblock {\em Convex Polytopes}, volume~221 of {\em Graduate Texts in Mathematics}.
\newblock Springer, New York, 2nd edition, 2003.
\newblock \href{https://doi.org/10.1007/978-1-4613-0019-9}{\path{doi:10.1007/978-1-4613-0019-9}}.

\bibitem{LP}
Arkadi Nemirovski.
\newblock {\em Introduction to Linear Optimization}.
\newblock World Scientific, Singapore, 2024.
\newblock Postprint dated March 17, 2024.
\newblock URL: \url{https://www2.isye.gatech.edu/~nemirovs/WSPrinted.pdf},
  \href {https://doi.org/10.1142/13457} {\path{doi:10.1142/13457}}.

\bibitem{RockafellarConvexAnalysis}
R.~Tyrrell Rockafellar.
\newblock {\em Convex Analysis}, volume~28 of {\em Princeton Mathematical
  Series}.
\newblock Princeton University Press, Princeton, NJ, 1970.

\bibitem{SalvettiSettepanella2007}
Mario Salvetti and Simona Settepanella.
\newblock Combinatorial {Morse} theory and minimality of hyperplane
  arrangements.
\newblock {\em Geometry \& Topology}, 11:1733--1766, 2007.
\newblock \href {https://doi.org/10.2140/gt.2007.11.1733}
  {\path{doi:10.2140/gt.2007.11.1733}}.

\end{thebibliography}
